\documentclass[11pt]{article}
\usepackage{enumerate}
\usepackage{amssymb,a4wide,latexsym,makeidx,epsfig,fleqn}
\usepackage{amsthm}
\usepackage{amsmath}
\usepackage{enumerate}
\usepackage{graphicx}
\usepackage{color}
\usepackage{epstopdf}
\newtheorem{theorem}{Theorem}[section]

\newtheorem{example}[theorem]{Example}
\newtheorem{definition}[theorem]{Definition}
\newtheorem{lemma}[theorem]{Lemma}

\begin{document}
\textwidth 150mm \textheight 225mm
\title{The row left rank of quaternion unit gain graphs in terms of pendant vertices   \footnote{This work is supported by the National Natural Science Foundations of China (No. 12371348, 12201258), the Postgraduate Research \& Practice Innovation Program of Jiangsu Normal University (No. 2024XKT1702).}}
\author{{ Yong Lu\footnote{Corresponding author.}, Qi Shen}\\
{\small  School of Mathematics and Statistics, Jiangsu Normal University,}\\ {\small  Xuzhou, Jiangsu 221116,
People's Republic
of China.}\\
{\small E-mails: luyong@jsnu.edu.cn, SQ\_jsnu@163.com }}
%%%%%%%%%%%%

\date{}
\maketitle
\begin{center}
\begin{minipage}{120mm}
\vskip 0.3cm
\begin{center}
{\small {\bf Abstract}}
\end{center}
{\small Let $\widetilde{G}=(G,U(\mathbb{Q}),\varphi)$  be a quaternion unit gain graph (or $U(\mathbb{Q})$-gain graph), where $G$ is the underlying graph of $\widetilde{G}$, $U(\mathbb{Q})=\{q\in \mathbb{Q}: |q|=1\}$  and $\varphi:\overrightarrow{E}\rightarrow U(\mathbb{Q})$ is the gain function such that $\varphi(e_{ij})=\varphi(e_{ji})^{-1}=\overline{\varphi(e_{ji})}$ for any adjacent vertices $v_{i}$ and $v_{j}$. Let $A(\widetilde{G})$ be the adjacency matrix of $\widetilde{G}$ and let $r(\widetilde{G})$ be the row left rank of $\widetilde{G}$. In this paper, we prove some lower bounds on the row left rank of $U(\mathbb{Q})$-gain graphs in terms of pendant vertices. All corresponding extremal graphs are characterized.

 }
\end{minipage}
\end{center}

\section{Introduction }
\emph{}In this paper, we consider only the graphs without multiedges and loops. Let $G=(V(G),E(G))$ be a simple graph, where $V(G)=\{v_{1},v_{2},\ldots,v_{n}\}$ and $E(G)$ are the vertex
set and the edge set of $G$, respectively. Let $v_{i},v_{j}\in V(G)$,  the \emph{adjacency matrix} $A(G)$ of $G$ is the symmetric $n\times n$ matrix with entries $a_{ij}=1$ if $v_{i}$ is adjacent to $v_{j}$ and $a_{ij}=0$ otherwise. The \emph{rank} (resp., \emph{nullity}) of $G$ is the rank (resp., nullity) of $A(G)$, denoted by $r(G)$ (resp., $\eta(G)$).

The research on the rank (or nullity) of graphs has attracted the attention of many scholars since it corresponds to the singularity of the graph. Collatz et al. \cite{CS} first proposed to characterize all graphs of order $n$ with $r(G)<n$.
Until today, this problem is still unsolved.
In the past decades, researchers have  focused on the boundaries of the nullities (or ranks) of graphs with given order in terms of
various graph parameters (and identifying the extremal graphs) such
as: the matching number  (see \cite{FHLL, LGUO, mf, RCZ, SST, WW}); the
number of pendant vertices (see
\cite{ccz, ctlz, MWTDAM,  wfg}); the maximum degree (see
\cite{clt, sl, WGT, WGUO, ZWS}); the girth (see \cite{candl, zwt}), etc.

Let $\mathbb{R}$ and $\mathbb{C}$ be the fields of the real numbers and complex numbers,
respectively. Let $\mathbb{Q}$ be a four-dimensional vector space over $\mathbb{R}$ with an
ordered basis, denoted by $1$, $i$, $j$, and $k$.
A \emph{real quaternion}, simply called \emph{quaternion}, is a vector
$q=x_{0}+x_{1}i+x_{2}j+x_{3}k\in \mathbb{Q},$
where $x_{0},x_{1},x_{2},x_{3}$ are real numbers and $i,j,k$ satisfy the following conditions:
%$$1q=q1=q, for~all~q\in \{i, j, k\};$$
$$i^{2}=j^{2}=k^{2}=-1;$$
$$ij=-ji=k, jk=-kj=i, ki=-ik=j.$$
 %If $a,b$ are any real numbers, while $\textbf{u},\textbf{v}$ are any two of $i,j,k$, then $(a\textbf{u})(b\textbf{v})=(ab)(\textbf{u}\textbf{v})$.

From \cite{ZF}, we know that if $x,y$ and $z$ are three different quaternions, then $(xy)^{-1}=y^{-1}x^{-1}$ and $(xy)z=x(yz)$. (Note that $xy\neq yx$, in general)

Let $q=x_{0}+x_{1}i+x_{2}j+x_{3}k\in \mathbb{Q}$. The \emph{conjugate} $\bar{q}$ (or $q^{\ast}$) of $q$ is $\bar{q}=x_{0}-x_{1}i-x_{2}j-x_{3}k$. The \emph{modulus} of $q$ is $|q|=\sqrt{q\bar{q}}=\sqrt{x_{0}^{2}+x_{1}^{2}+x_{2}^{2}+x_{3}^{2}}$. If $q\neq 0$, then the \emph{inverse} of $q$ is $q^{-1}=\frac{\bar{q}}{|q|^{2}}$. The \emph{real part} of $q$ is $Re(q)=x_{0}$. The \emph{imaginary part} of $q$ is $Im(q)=x_{1}i+x_{2}j+x_{3}k$. The \emph{row left (right) rank} of a quaternion matrix $A\in \mathbb{Q}^{m\times n}$ is the maximum number of rows of $A$ that are left (right) linearly independent. The \emph{column left (right) rank} of a quaternion matrix $A\in \mathbb{Q}^{m\times n}$ is the maximum number of columns of $A$ that are left (right) linearly independent.

An oriented edge from $v_{i}$ to $v_{j}$ is denoted by $e_{ij}$.
Thus $e_{ij}$ and $e_{ji}$ are considered to be distinct.
Let $\overrightarrow{E}$ denote the set of $\{e_{ij}, e_{ji}:~v_{i}v_{j}\in E\}$. Even though $e_{ij}$ stands for an edge and an oriented edge simultaneously, it will always be clear from the context.
A \emph{gain graph} is a graph with the following additional structure: each orientation of an edge is given a group element, called a \emph{gain}, which is the inverse of the group element assigned to the opposite orientation.
Belardo et al. \cite{BBCR} studied quaternion unit gain graphs and their associated spectral theories.
Denote by $\widetilde{G}=(G,U(\mathbb{Q}),\varphi)$  a quaternion unit gain graph (or $U(\mathbb{Q})$-gain graph), where $G$ is the \emph{underlying graph} of $\widetilde{G}$, $U(\mathbb{Q})=\{q\in \mathbb{Q}: |q|=1\}$ and $\varphi:\overrightarrow{E}\rightarrow U(\mathbb{Q})$ is the gain function such that $\varphi(e_{ij})=\varphi(e_{ji})^{-1}=\overline{\varphi(e_{ji})}$.
For convenience, $\varphi(e_{ij})$ is also written as $\varphi_{v_{i}v_{j}}$ for $e_{ij}\in E(\widetilde{G})$.
The \emph{adjacency matrix} of  $\widetilde{G}$ is the  Hermitian matrix $A(\widetilde{G})=(h_{ij})_{n\times n}$, where $h_{ij}=\varphi(e_{ij})=\varphi_{v_{i}v_{j}}$ if $e_{ij}\in E(\widetilde{G})$, and $h_{ij}=0$ otherwise.
The row left rank $r(\widetilde{G})$ of $A(\widetilde{G})$ is called the \emph{rank}  of $\widetilde{G}$.

 Note that quaternions do not satisfy the commutative law of multiplication. By the following lemma and  the example in \cite{QNZ}, we know that the row left rank of a
quaternion matrix is not necessarily equal to the row right rank (the same as column left rank and column right rank).

\noindent\begin{lemma}\label{lem:2.1}\cite{QNZ}
The row left rank of a quaternion matrix $A$ equals the column right rank of $A$. The row right rank of a quaternion matrix $A$ equals the column left rank of $A$.
\end{lemma}

Note that quaternion unit gain graphs will generalize simple graphs ($\varphi(\overrightarrow{E})=\{1\}$),
signed graphs ($\varphi(\overrightarrow{E})\subseteq\{1,-1\}$),
mixed graphs ($\varphi(\overrightarrow{E})\subseteq\{1,i,-i\}$)
and complex unit gain graphs ($\varphi(\overrightarrow{E})\subseteq\{z\in \mathbb{C}: |z|=1\}$).
Recently, people have extended the research of the rank (or nullity) of simple graphs to
signed graphs (see \cite{BB, HHL, LWZ, lwn, WS, wlt})
and complex unit gain graphs (see \cite{hhd, hhy, landy, LUWH, lwx, LWZ1, landwu, REFF, xzw, yqt}).
Zhou and Lu \cite{QNZ} obtained the relationship between the row left rank of a quaternion unit gain graph and the rank of its underlying graph.

For a vertex $x\in V(\widetilde{G})$, we define $N_{\widetilde{G}}(x)=\{y\in V(\widetilde{G}): e_{xy},e_{yx}\in E(\widetilde{G})\}$, which is called the \emph{neighbor set} of $x$ in $\widetilde{G}$. The \emph{degree} of $x$ is the number of vertices which are adjacent to $x$, denoted by $d_{\widetilde{G}}(x)$. If $d_{\widetilde{G}}(x)=1$, then $x$ is called a \emph{pendant vertex} (or a \emph{leaf}). The number of leaves in $\widetilde{G}$ is denoted by $p(\widetilde{G})$. If $p(\widetilde{G})=0$, then we call $\widetilde{G}$ a \emph{leaf-free graph}. If $d_{\widetilde{G}}(x)\geq2$, then $x$ is called a \emph{major vertex}. Let $c(\widetilde{G})$ be the \emph{cyclomatic number} of a $U(\mathbb{Q})$-gain graph $\widetilde{G}$, that is $c(\widetilde{G})=|E(\widetilde{G})|-|V(\widetilde{G})|+\omega(\widetilde{G})$, where $\omega(\widetilde{G})$ is the number of connected components of $\widetilde{G}$.
Denote by $\widetilde{P}_{n}$ and $\widetilde{C}_{n}$ a $U(\mathbb{Q})$-gain path and a $U(\mathbb{Q})$-gain cycle on $n$ vertices, respectively.
For $u,v\in V(\widetilde{G})$, we denote by $d_{\widetilde{G}}(u,v)$ the length of a shortest path between $u$ and $v$ in $\widetilde{G}$.

In this paper, our main results consist in getting some lower bounds of the row left rank of $U(\mathbb{Q})$-gain graphs in terms of $c(\widetilde{G})$ and $p(\widetilde{G})$ (see the following theorem), and characterizing the corresponding extremal $U(\mathbb{Q})$-gain graphs (see Theorems \ref{th:5.1.}, \ref{th:5.3.} and \ref{th:5.5.}). The results of this paper are a generalization of those in \cite{MWTDAM}.

\noindent\begin{theorem}\label{th:3.2.}
Let $\widetilde{G}$ be a $U(\mathbb{Q})$-gain graph of order $n$ in which every component of $\widetilde{G}$ has at least two vertices. Then
\begin{align*}
r(\widetilde{G})\geq\left\{\begin{array}{ll}
n-2c(\widetilde{G})-p(\widetilde{G})+1,  ~~if~p(\widetilde{G})\geq1;\\
n-2c(\widetilde{G}), ~~if\;p(\widetilde{G})=0\;and\;\widetilde{G}\;is\;a\;cycle \text{-} disjoint\;U(\mathbb{Q})\text{-}gain\;graph;\\
n-2c(\widetilde{G})+1, ~~if\;p(\widetilde{G})=0\;and\;some\;cycles\;have\;common\;vertices.\\
\end{array}\right.
\end{align*}
\end{theorem}

Let $\widetilde{G}$ be a $U(\mathbb{Q})$-gain graph with vertex set $V(\widetilde{G})$ and let $U\subseteq V(\widetilde{G})$, denote by $\widetilde{G}-U$ the $U(\mathbb{Q})$-gain graph obtained from $\widetilde{G}$ by removing the vertices in $U$ together with all incident edges. When $U=\{x\}$, we write $\widetilde{G}-U$ simply as $\widetilde{G}-x$.
A vertex $x\in V(\widetilde{G})$ is called a \emph{cut-point} of a $U(\mathbb{Q})$-gain graph $\widetilde{G}$ if $\widetilde{G}-x$ has more connected components than $\widetilde{G}$.
A \emph{block} of a $U(\mathbb{Q})$-gain graph is a maximal connected graph which does not have any cut-point.
Sometimes we use the notation $\widetilde{G}-\widetilde{H}$ instead of $\widetilde{G}-V(\widetilde{H})$ when $\widetilde{H}$ is an induced subgraph of $\widetilde{G}$.
If $\widetilde{G}_{1}$ is an induced subgraph of $\widetilde{G}$ and $x$ is a vertex not in $\widetilde{G}_{1}$, we write the subgraph of $\widetilde{G}$ induced by $V(\widetilde{G}_{1})\cup\{x\}$ simply as $\widetilde{G}_{1}+x$.

\section{Preliminaries}

By basic matrix theory, we have the following lemma.

\noindent\begin{lemma}\label{le:2.11.}\cite{QNZ}
Let $\widetilde{G}=\widetilde{G}_{1}\cup \widetilde{G}_{2}\cup \cdots \cup \widetilde{G}_{t}$, where $\widetilde{G}_{1},\widetilde{G}_{2},\cdots,\widetilde{G}_{t}$ are connected components of  $\widetilde{G}$. Then $r(\widetilde{G})=\sum^{t}_{i=1}r(\widetilde{G}_{i})$.
\end{lemma}

Assume $\widetilde{G}$ is connected. Let $\widetilde{T}$ be a spanning tree and let $u$ be a root vertex.
For $v\in V(\widetilde{G})$, let $\widetilde{P}_{vu}$ be the unique path in $\widetilde{T}$ from $v$ to $u$.
Let $\theta: V\rightarrow U(\mathbb{Q})$ be a switching function such that $\theta(v)=\varphi(\widetilde{P}_{vu})$.
Switching the $U(\mathbb{Q})$-gain graph $\widetilde{G}$ by $\theta$ means forming a new quaternion unit gain graph $\widetilde{G}^{\theta}$, whose underlying graph is the same as $\widetilde{G}$, but whose gain function is defined on an edge $xy$ by $\varphi^{\theta}(xy)= \theta(x)^{-1}\varphi(xy)\theta(y)$.
If there exists a switching function $\theta$ such that $\widetilde{G}_{2}=\widetilde{G}_{1}^{\theta}$, then $\widetilde{G}_{1}$ and $\widetilde{G}_{2}$ are called \emph{switching equivalent}, denoted by $\widetilde{G}_{1}\leftrightarrow \widetilde{G}_{2}$.
Note that two switching equivalent $U(\mathbb{Q})$-gain graphs have the same rank.

\noindent\begin{lemma}\label{le:2.0.}\cite{QNZ}
Let $\widetilde{T}$ be a $U(\mathbb{Q})$-gain tree of order $n$. Then $A(\widetilde{T})$ and $A(T)$ have the same rank.
\end{lemma}

According to the above lemma, we get that the spectral theory of $U(\mathbb{Q})$-gain trees is indistinguishable from the spectral theory of their underlying graphs. In fact, the gain map $\varphi$ does not play any role in spectral results on trees.

\noindent\begin{lemma}\label{le:2.1.}\cite{QNZ}
Let $\widetilde{P}$ be a $U(\mathbb{Q})$-gain path of order $n$. If $n$ is odd, then $r(\widetilde{P})=n-1$; if $n$ is even, then $r(\widetilde{P})=n$.
\end{lemma}

\noindent\begin{definition}\label{de:2.2.}\cite{QNZ}
Let $\widetilde{C}_{n} (n\geq3)$ be a $U(\mathbb{Q})$-gain cycle and let $$\varphi(\widetilde{C}_{n})=\varphi_{v_{1}v_{2}}\varphi_{v_{2}v_{3}}\cdots \varphi_{v_{n-1}v_{n}}\varphi_{v_{n}v_{1}}.$$ Then $\widetilde{C}_{n}$ is said to be:
\begin{displaymath}
\left\{\
        \begin{array}{ll}
          \rm Type~1,&  \emph{if}~\varphi(\widetilde{C}_{n})=(-1)^{n/2}~\emph{and}~n~\emph{is~even};\\
          \rm Type~2,& \emph{if}~\varphi(\widetilde{C}_{n})\neq(-1)^{n/2}~\emph{and}~n~\emph{is~even};\\
          \rm Type~3,& \emph{if}~Re\left((-1)^{{(n-1)}/{2}}\varphi(\widetilde{C}_{n})\right)\neq 0~\emph{and}~n~\emph{is~odd};\\
          \rm Type~4,& \emph{if}~Re\left((-1)^{{(n-1)}/{2}}\varphi(\widetilde{C}_{n})\right)=0~\emph{and}~n~\emph{is~odd}.
        \end{array}
      \right.
\end{displaymath}
\end{definition}

\noindent\begin{lemma}\label{le:2.3.}\cite{QNZ}
Let $\widetilde{C}_{n}$ be a $U(\mathbb{Q})$-gain cycle of order $n$. Then
\begin{align*}
r(\widetilde{C}_{n})=\left\{\begin{array}{ll}
n-2, & if \;\widetilde{C}_{n}\; is \;of\; Type \;1;\\
n,& if \;\widetilde{C}_{n}\; is \;of\; Type \;2 \;or\; 3;\\
n-1,& if \;\widetilde{C}_{n}\; is \;of\; Type \;4.\\
\end{array}\right.
\end{align*}
\end{lemma}

\noindent\begin{lemma}\label{le:2.10.}\cite{QNZ}
Let $\widetilde{G}$ be a $U(\mathbb{Q})$-gain graph obtained by identifying a vertex of a $U(\mathbb{Q})$-gain cycle $\widetilde{C}_{n}$ with a vertex of a $U(\mathbb{Q})$-gain graph $\widetilde{G}_{1}$ of order $m~(m\geq 1)$ $(V(\widetilde{C}_{n})\cap V(\widetilde{G}_{1})=\{u\})$ and let $\widetilde{G}_{2}=\widetilde{G}_{1}-u$. Then
\begin{align*}
\left\{\begin{array}{ll}
r(\widetilde{G})=n-2+r(\widetilde{G}_{1}), & if \;\widetilde{C}_{n}\; is \;of\; Type \;1;\\
r(\widetilde{G})=n+r(\widetilde{G}_{2}),& if \;\widetilde{C}_{n}\; is \;of\; Type \;2;\\
r(\widetilde{G})=n-1+r(\widetilde{G}_{1}),& if \;\widetilde{C}_{n}\; is \;of\; Type \;4;\\
n-1+r(\widetilde{G}_{2})\leq r(\widetilde{G})\leq n+r(\widetilde{G}_{1}),& if \;\widetilde{C}_{n}\; is \;of\; Type \;3.\\
\end{array}\right.
\end{align*}
\end{lemma}

The remaining lemmas below are the result of the effect of vertex and edge operations on the  row left rank of  quaternion unit gain graphs.

\noindent\begin{lemma}\label{le:2.4.}\cite{QNZ}
Let $\widetilde{G}$ be a $U(\mathbb{Q})$-gain graph and let $v$ be a vertex of $\widetilde{G}$. Then $r(\widetilde{G})-2\leq r(\widetilde{G}-v)\leq r(\widetilde{G})$.
\end{lemma}

\noindent\begin{lemma}\label{le:2.5.}\cite{QNZ}
Let $\widetilde{G}$ be a $U(\mathbb{Q})$-gain graph and let $x$ be a pendant vertex of $\widetilde{G}$. If $y$ is adjacent to  $x$ in $\widetilde{G}$, then $r(\widetilde{G})=r(\widetilde{G}-x-y)+2$.
\end{lemma}

\noindent\begin{lemma}\label{le:4.1.}
Let $\widetilde{H}$ and $\widetilde{K}$ be two $U(\mathbb{Q})$-gain graphs. If $\widetilde{G}$ is obtained by identifying a vertex $v$ of $\widetilde{H}$ with a vertex $u$ of $\widetilde{K}$, then $r(\widetilde{G})\geq r(\widetilde{K})+r(\widetilde{H}-v)-1$.
\end{lemma}
\noindent\textbf{Proof.}
Let $V(\widetilde{H})=\{v_{1},v_{2},\ldots,v_{n}\}$ and let $V(\widetilde{K})=\{u_{1},u_{2},\ldots,u_{m}\}$.
Without loss of generality, let $v_{n}=v$ and let $u_{1}=u$.
Let $h_{ij}=\varphi'_{v_{i}v_{j}}$ $(1\leq i, j\leq n)$ and let $p_{ij}=\varphi''_{u_{i}u_{j}}$ $(1\leq i, j\leq m)$, where $\varphi'$ and $\varphi''$ are the gain functions defining $\widetilde{H}$ and $\widetilde{K}$ respectively.
\begin{align*}
 A(\widetilde{G})= \left (
\begin{array}{cccccccccc}
 A(\widetilde{H}-v) & \alpha & \textbf{0}_{(n-1)\times(m-1)}\\
 \overline{\alpha}^{T} & 0 & \beta\\
 \textbf{0}_{(m-1)\times(n-1)} & \overline{\beta}^{T} & A(\widetilde{K}-u)\\
 \end{array}
 \right),
 \end{align*}
where $\alpha=(h_{1n},h_{2n},\ldots,h_{(n-1)n})^{T}$ and $\beta=(p_{12},p_{13},\ldots,p_{1m})$.
Since
\begin{align*}
 A(\widetilde{G})+\left (
\begin{array}{cccccccccc}
 \textbf{0}_{(n-1)\times(n-1)} & \textbf{0}_{(n-1)\times1} & \textbf{0}_{(n-1)\times(m-1)}\\
 -\overline{\alpha}^{T} & 0 & \textbf{0}_{1\times(m-1)}\\
 \textbf{0}_{(m-1)\times(n-1)} & \textbf{0}_{(m-1)\times1} & \textbf{0}_{(m-1)\times(m-1)}\\
 \end{array}
 \right)=\left (
\begin{array}{cccccccccc}
  A(\widetilde{H}-v) & \alpha & \textbf{0}_{(n-1)\times(m-1)}\\
 \textbf{0}_{1\times(n-1)} & 0 & \beta\\
 \textbf{0}_{(m-1)\times(n-1)} & \overline{\beta}^{T} & A(\widetilde{K}-u)\\
 \end{array}
 \right),
 \end{align*}
the row left rank of the left side matrix is at most $r(\widetilde{G})+1$ and the row left rank of the right side matrix is at least $r(\widetilde{K})+r(\widetilde{H}-v)$, then $r(\widetilde{G})\geq r(\widetilde{K})+r(\widetilde{H}-v)-1$.
\quad $\square$~\\

\noindent\begin{lemma}\label{le:4.2.}
Let $\widetilde{H}$ and $\widetilde{K}$ be two $U(\mathbb{Q})$-gain graphs. If $\widetilde{G}$ is obtained from $\widetilde{H}$ and $\widetilde{K}$ by connecting one vertex $v_{1}$ of $\widetilde{H}$ and one vertex $v_{t}$ of $\widetilde{K}$ with a $U(\mathbb{Q})$-gain path $\widetilde{P}_{t}~(t\geq2)$, then $r(\widetilde{G})\geq r(\widetilde{H})+r(\widetilde{K})+t-3$.
\end{lemma}
\noindent\textbf{Proof.}
If $t=2$, then $\widetilde{G}$ is obtained by identifying a vertex of $\widetilde{H}$ with a vertex of $\widetilde{K}+v_{1}$ ($V(\widetilde{H})\cap V(\widetilde{K}+v_{1})=\{v_{1}\}$).
By Lemma \ref{le:4.1.}, we have $r(\widetilde{G})\geq r(\widetilde{H})+r(\widetilde{K}+v_{1}-v_{1})-1=r(\widetilde{H})+r(\widetilde{K})-1=r(\widetilde{H})+r(\widetilde{K})+t-3$.

If $t=3$, then $V(\widetilde{P}_{t})=V(\widetilde{P}_{3})=\{v_{1},v_{2},v_{3}\}$. By Lemma \ref{le:2.4.},  $r(\widetilde{G})\geq r(\widetilde{G}-v_{2})=r(\widetilde{H})+r(\widetilde{K})=r(\widetilde{H})+r(\widetilde{K})+t-3$.

Suppose that $t\geq4$, let $V(\widetilde{P}_{t})=\{v_{1},v_{2},\ldots,v_{t}\}$ ($v_{1}\in V(\widetilde{H})$ and $v_{t}\in V(\widetilde{K})$) and let $h_{i}=\varphi_{v_{i}v_{i+1}}$ $(1\leq i\leq t-1)$.
\begin{align*}
 A(\widetilde{G})= \left (
\begin{array}{ccccccccccccc}
 A(\widetilde{H}) & \alpha & \textbf{0}_{m\times1} & \cdots & \textbf{0}_{m\times1} & \textbf{0}_{m\times1} & \textbf{0}_{m\times n}\\
 \overline{\alpha}^{T} & 0 & h_{2} & \cdots & 0 & 0 & \textbf{0}_{1\times n}\\
 \textbf{0}_{1\times m} & \overline{h_{2}} & 0 & \cdots & 0 & 0 & \textbf{0}_{1\times n}\\
 \vdots & \vdots & \vdots &\ddots & \vdots & \vdots & \vdots\\
 \textbf{0}_{1\times m} & 0 & 0 & \cdots & 0 & h_{t-2} & \textbf{0}_{1\times n}\\
 \textbf{0}_{1\times m} & 0 & 0 & \cdots & \overline{h_{t-2}} & 0 & \beta\\
 \textbf{0}_{n\times m} & \textbf{0}_{n\times1} & \textbf{0}_{n\times1} & \cdots & \textbf{0}_{n\times1} & \overline{\beta}^{T} & A(\widetilde{K})\\
 \end{array}
 \right),
 \end{align*}
where $|V(\widetilde{H})|=m$, $|V(\widetilde{K})|=n$, $\alpha=(0,\ldots,0,h_{1})^{T}$ and $\beta=(h_{t-1},0,\ldots,0)$.

Let $B$ be the truncated matrix of block matrix $A(\widetilde{G})$ obtained by deleting the second row $(\overline{\alpha}^{T},0,h_{2},\cdots,0,0, \textbf{0}_{1\times n})$ and the second last column $(\textbf{0}_{m\times1},0,0,\cdots,h_{t-2},0, \overline{\beta}^{T})^{T}$, that is
\begin{align*}
 B= \left (
\begin{array}{ccccccccccccc}
 A(\widetilde{H}) & \alpha & \textbf{0}_{m\times1} & \cdots & \textbf{0}_{m\times1} & \textbf{0}_{m\times n}\\
 \textbf{0}_{1\times m} & \overline{h_{2}} & 0 & \cdots & 0 & \textbf{0}_{1\times n}\\
 \vdots & \vdots & \vdots &\ddots & \vdots & \vdots\\
 \textbf{0}_{1\times m}& 0 & 0 & \cdots & 0 & \textbf{0}_{1\times n}\\
 \textbf{0}_{1\times m} & 0 & 0 & \cdots & \overline{h_{t-2}} & \beta\\
 \textbf{0}_{n\times m} & \textbf{0}_{n\times1} & \textbf{0}_{n\times1} & \cdots & \textbf{0}_{n\times1} & A(\widetilde{K})\\
 \end{array}
 \right).
 \end{align*}
Since $B$ is a block upper triangular matrix, we have
$$r(\widetilde{G})\geq r(B)\geq r(\widetilde{H})+r(\widetilde{K})+t-3.~~\square$$

Without compromising the outcome, denote by $\widetilde{G}-e_{ij}$ the $U(\mathbb{Q})$-gain graph obtained from $\widetilde{G}$ by removing the edge $e_{ij}$ or $e_{ji}$ (no matter $e_{ij}$ or $e_{ji}$ $\in E(\widetilde{G})$).

\noindent\begin{lemma}\label{le:2.6.}
 Let $\widetilde{G}$ be a $U(\mathbb{Q})$-gain graph and let $e_{ij}$ be an edge of $\widetilde{G}$. Then $r(\widetilde{G})\geq r(\widetilde{G}-e_{ij})-2$.
\end{lemma}
\noindent\textbf{Proof.}
Let $V(\widetilde{G})=\{v_{1},v_{2},\ldots,v_{n}\}$ and let $v_{1},v_{2}$ be adjacent.
Without loss of generality, let $e_{ij}=e_{12}$ and let $h_{ij}=\varphi_{v_{i}v_{j}}$ $(1\leq i, j\leq n)$.
\begin{align*}
 A(\widetilde{G})= \left (
\begin{array}{ccccccccccccc}
 0 & h_{12} & \alpha\\
 h_{21} & 0 & \beta\\
 \overline{\alpha}^{T} & \overline{\beta}^{T} & A(\widetilde{G}-v_{1}-v_{2})\\
 \end{array}
 \right),
 \end{align*}
where $\alpha=(h_{13},h_{14},\ldots,h_{1n})$ and $\beta=(h_{23},h_{24},\ldots,h_{2n})$. Then
\begin{align*}
A(\widetilde{G}-e_{12})=\left (
\begin{array}{ccccccccccccc}
 0 & -h_{12} & \textbf{0}_{1\times (n-2)}\\
 -h_{21} & 0 & \textbf{0}_{1\times (n-2)}\\
 \textbf{0}_{(n-2)\times 1} & \textbf{0}_{(n-2)\times 1} & \textbf{0}_{(n-2)\times (n-2)}\\
 \end{array}
 \right)+A(\widetilde{G})
=\left (
\begin{array}{ccccccccccccc}
 0 & 0 & \alpha\\
 0 & 0 & \beta\\
 \overline{\alpha}^{T} & \overline{\beta}^{T} & A(\widetilde{G}-v_{1}-v_{2})\\
 \end{array}
 \right).
 \end{align*}
Thus
$r(\widetilde{G})\geq r(\widetilde{G}-e_{12})-2.~~\square$

\noindent\begin{lemma}\label{le:2.7.}
Suppose $\widetilde{G}$ be a $U(\mathbb{Q})$-gain graph and $\widetilde{P}_{6}=v_{1}v_{2}v_{3}v_{4}v_{5}v_{6}$ be a $U(\mathbb{Q})$-gain path with $d_{\widetilde{G}}(v_{i})=2$ $(i=2,3,4,5)$ in $\widetilde{G}$.
Let $\widetilde{G}_{1}$ be a $U(\mathbb{Q})$-gain graph obtained by replacing $\widetilde{P}_{6}$ with a new edge $v_{1}v_{6}$ with gain $\varphi(\widetilde{P}_{6})$. Then $r(\widetilde{G})=r(\widetilde{G}_{1})+4$.
\end{lemma}
\noindent\textbf{Proof.}
Let $\widetilde{H}=\widetilde{G}-\widetilde{P}_{6}$ and let $V(\widetilde{H})=\{u_{1},u_{2},\ldots,u_{n-6}\}$.
Let $h_{i}=\varphi_{v_{i}v_{i+1}}$ $(1\leq i\leq5)$ and let $q_{ij}=\varphi_{u_{i}v_{j}}$ $(1\leq i\leq n-6~and~j=1,6)$.
\begin{align*}
 A(\widetilde{G})=\left (
\begin{array}{ccccccccccccc}
 A(\widetilde{H}) & \alpha & \textbf{0}_{(n-6)\times1} & \textbf{0}_{(n-6)\times1} & \textbf{0}_{(n-6)\times1} & \textbf{0}_{(n-6)\times1} & \beta\\
 \overline{\alpha}^{T} & 0 & h_{1} & 0 & 0 & 0 & 0\\
 \textbf{0}_{1\times(n-6)} & \overline{h}_{1} & 0 & h_{2} & 0 & 0 & 0 \\
 \textbf{0}_{1\times(n-6)} & 0 & \overline{h}_{2} & 0 & h_{3} & 0 & 0 \\
 \textbf{0}_{1\times(n-6)} & 0 & 0 & \overline{h}_{3} & 0 & h_{4} & 0 \\
 \textbf{0}_{1\times(n-6)} & 0 & 0 & 0 & \overline{h}_{4} & 0 & h_{5}\\
  \overline{\beta}^{T} & 0 & 0 & 0 & 0 & \overline{h}_{5} & 0\\
 \end{array}
 \right),
 \end{align*}
where $\alpha=(q_{11},q_{21},\ldots,q_{(n-6)1})^{T}$ and $\beta=(q_{16},q_{26},\ldots,q_{(n-6)6})^{T}$.
We multiply $-h_{1}h_{2}$ on the left side of 4-th row and add it to 2-th row,
and multiply $-\overline{h}_{5}\overline{h}_{4}$ on the left side of 5-th row and add it to 7-th row:
\begin{align*}
 r(\widetilde{G})=r\left (
 \begin{array}{ccccccc}
  A(\widetilde{H}) & \alpha & \textbf{0}_{(n-6)\times 1} & \textbf{0}_{(n-6)\times 1} & \textbf{0}_{(n-6)\times 1} & \textbf{0}_{(n-6)\times 1} & \beta\\
 \overline{\alpha}^{T} & 0 & 0 & 0 & -h_{1}h_{2}h_{3} & 0 & 0\\
 \textbf{0}_{1\times(n-6)} & \overline{h}_{1} & 0 & h_{2} & 0 & 0 & 0 \\
 \textbf{0}_{1\times(n-6)} & 0 & \overline{h}_{2} & 0 & h_{3} & 0 & 0 \\
 \textbf{0}_{1\times(n-6)} & 0 & 0 & \overline{h}_{3} & 0 & h_{4} & 0 \\
 \textbf{0}_{1\times(n-6)} & 0 & 0 & 0 & \overline{h}_{4} & 0 & h_{5}\\
  \overline{\beta}^{T} & 0 & 0 & -\overline{h}_{5}\overline{h}_{4}\overline{h}_{3} & 0 & 0 & 0\\
 \end{array}
 \right).
 \end{align*}
We multiply $h_{1}h_{2}h_{3}h_{4}$ on the left side of 6-th row and add it to 2-th row, and multiply $\overline{h}_{5}\overline{h}_{4}\overline{h}_{3}\overline{h}_{2}$ on the left side of 3-th row and add it to 7-th row:
\begin{align*}
 r(\widetilde{G})=r\left (
 \begin{array}{ccccccc}
  A(\widetilde{H}) & \alpha & \textbf{0}_{(n-6)\times 1} & \textbf{0}_{(n-6)\times 1} & \textbf{0}_{(n-6)\times 1} & \textbf{0}_{(n-6)\times 1} & \beta\\
 \overline{\alpha}^{T} & 0 & 0 & 0 & 0 & 0 & h_{1}h_{2}h_{3}h_{4}h_{5}\\
 \textbf{0}_{1\times(n-6)} & \overline{h}_{1} & 0 & h_{2} & 0 & 0 & 0 \\
 \textbf{0}_{1\times(n-6)} & 0 & \overline{h}_{2} & 0 & h_{3} & 0 & 0 \\
 \textbf{0}_{1\times(n-6)} & 0 & 0 & \overline{h}_{3} & 0 & h_{4} & 0 \\
 \textbf{0}_{1\times(n-6)} & 0 & 0 & 0 & \overline{h}_{4} & 0 & h_{5}\\
  \overline{\beta}^{T} & \overline{h}_{5}\overline{h}_{4}\overline{h}_{3}\overline{h}_{2}\overline{h}_{1} & 0 & 0 & 0 & 0 & 0\\
 \end{array}
 \right).
 \end{align*}
Let $h'_{16}=h_{1}h_{2}h_{3}h_{4}h_{5}=\varphi(\widetilde{P}_{6})$. Write
\begin{align*}
 B= \left (
\begin{array}{ccccccccccccc}
 A(\widetilde{H}) & \alpha & \textbf{0}_{(n-6)\times 1} & \textbf{0}_{(n-6)\times 1} & \textbf{0}_{(n-6)\times 1} & \textbf{0}_{(n-6)\times 1} & \beta\\
 \overline{\alpha}^{T} & 0 & 0 & 0 & 0 & 0 & h'_{16}\\
  \overline{\beta}^{T} & \overline{h}'_{16} & 0 & 0 & 0 & 0 & 0 \\
 \end{array}
 \right),
 \end{align*}
 \begin{align*}
 A(\widetilde{G}_{1})=\left (
\begin{array}{ccccccccccccc}
 A(\widetilde{H}) & \alpha & \beta\\
 \overline{\alpha}^{T} & 0 & h'_{16}\\
  \overline{\beta}^{T} & \overline{h}'_{16} & 0 \\
 \end{array}
 \right).
 \end{align*}
Then $r(\widetilde{G})=r(B)+4$ and $r(B)=r(\widetilde{G}_{1})$. Thus $r(\widetilde{G})=r(\widetilde{G}_{1})+4$.
$\square$

\section{Lower bounds on the row left rank of $\widetilde{G}$}
In this section, we will obtain some lower bounds of the row left rank of $U(\mathbb{Q})$-gain graphs
in terms of the cyclomatic number and the number of pendant vertices.

Let $x$ be a vertex of a $U(\mathbb{Q})$-gain graph $\widetilde{G}$. If $y$ is adjacent to $x$ and $d_{\widetilde{G}}(y)=2$, then $y$ is called a \emph{2-degree neighbor} of $x$. By using the similar arguments of the proof of Lemma 2.6 in \cite{MWTDAM}, we get the following result.
%%%%%%%%%%%%%%
\noindent\begin{lemma}\label{le:3.1.}
Let $\widetilde{G}$ be a $U(\mathbb{Q})$-gain graph with a vertex $x$.
Let $r$ be the number of components containing 2-degree neighbors of $x$,
$m$ be the number of 2-degree neighbors of $x$,
and let $s$ be the number of components of $\widetilde{G}-x$. We have
\begin{enumerate}[(1)]
\item $d_{\widetilde{G}}(x)+r\geq m+s$;
\item If x lies on a $U(\mathbb{Q})$-gain cycle of $\widetilde{G}$, then $2d_{\widetilde{G}}(x)+r\geq m+2s+1$;
\item $c(\widetilde{G}-x)=c(\widetilde{G})-d_{\widetilde{G}}(x)+s$.
\end{enumerate}
\end{lemma}

Now, we get some lower bounds of the row left rank $r(\widetilde{G})$ of $U(\mathbb{Q})$-gain graphs.

\noindent\textbf{Proof of Theorem \ref{th:3.2.}.}
We may assume that $\widetilde{G}$ is connected. We proceed by induction on the order of $\widetilde{G}$. If $n=2$, then $r(\widetilde{G})=2$, $c(\widetilde{G})=0$ and $p(\widetilde{G})=2$, thus $r(\widetilde{G})> n-2c(\widetilde{G})-p(\widetilde{G})+1$. Now we assume that $n\geq3$ and the inequality holds for any connected $U(\mathbb{Q})$-gain graph $\widetilde{G}$ with order $2\leq|V(\widetilde{G})|\leq n-1$. Now we discuss two cases.

\textbf{Case 1.} $p(\widetilde{G})=0$.

In this case, $c(\widetilde{G})\geq 1$. Let $x$ be a vertex lying on a $U(\mathbb{Q})$-gain cycle
and let $\widetilde{H}_{1},\widetilde{H}_{2},\ldots, \widetilde{H}_{s}$ be the connected components of $\widetilde{G}-x$.
Without loss of generality, let $\widetilde{H}_{i}$ be the components containing 2-degree neighbors of $x$ for $i=1,...,r$,
and let $\widetilde{H}_{j}$ be the components containing  no 2-degree neighbors of $x$ for $j=r+1,...,s$.
This arrangement implies that $\widetilde{H}_{i}$ has pendant vertices for $i=1,...,r$ and
$\widetilde{H}_{j}$ has no pendant vertices for $j=r+1,...,s$.
Since $\widetilde{G}$ has no pendant vertices, $2\leqslant|V(\widetilde{H}_{i})|<|V(\widetilde{G})|$ for $i=1,...,s$, by the induction hypothesis, we get
$$r(\widetilde{H}_{i})\geq |V(\widetilde{H}_{i})|-2c(\widetilde{H}_{i})-p(\widetilde{H}_{i})+1~for~i=1,...,r;$$
$$r(\widetilde{H}_{j})\geq |V(\widetilde{H}_{j})|-2c(\widetilde{H}_{j})~for~j=r+1,...,s.$$
By Lemmas \ref{le:2.4.} and \ref{le:2.11.}, we obtain
\begin{align}
\nonumber
r(\widetilde{G})&\geq r(\widetilde{G}-x)\nonumber
\\&=\sum^{r}_{i=1}r(\widetilde{H}_{i})+\sum^{s}_{j=r+1}r(\widetilde{H}_{j})\nonumber
\\&\geq\sum^{r}_{i=1}(|V(\widetilde{H}_{i})|-2c(\widetilde{H}_{i})-p(\widetilde{H}_{i})+1)+\sum^{s}_{j=r+1}(|V(\widetilde{H}_{j})|-2c(\widetilde{H}_{i}))\nonumber
\\&=n-1+r-2\sum^{s}_{i=1}c(\widetilde{H}_{i})-\sum^{r}_{i=1}p(\widetilde{H}_{i})\nonumber
\\&=n-1+r-2c(\widetilde{G}-x)-p(\widetilde{G}-x).\nonumber
\end{align}
Since $p(\widetilde{G})=0$, $p(\widetilde{G}-x)=m$, where $m$ denotes the number of 2-degree neighbors of $x$. By Lemma \ref{le:3.1.}(3),
\begin{equation}
r(\widetilde{G})\geq n-2c(\widetilde{G})+[2d_{\widetilde{G}}(x)+r-m-2s-1].
\end{equation}
Combining this with Lemma \ref{le:3.1.}(2), we obtain $r(\widetilde{G})\geq n-2c(\widetilde{G})$.
Next, we discuss the subcases based on whether the cycles in $\widetilde{G}$ have common vertices.

\textbf{Subcase 1.1.} Some $U(\mathbb{Q})$-gain cycles of $\widetilde{G}$ have common vertices.
Let $\widetilde{C}_{1}$ and $\widetilde{C}_{2}$ be two distinct $U(\mathbb{Q})$-gain cycles of $\widetilde{G}$ have common vertices,
and let $x$ be a common vertex of $\widetilde{C}_{1}$ and $\widetilde{C}_{2}$ such that $N_{\widetilde{C}_{1}}(x)\neq N_{\widetilde{C}_{2}}(x)$.
Hence $d_{\widetilde{G}}(x)\geq s+2$. By  Lemma \ref{le:3.1.}(1), we obtain
\begin{equation}
2d_{\widetilde{G}}(x)+r\geqslant 2s+m+2.
\end{equation}
Combining (1) and (2), we obtain $r(\widetilde{G})\geqslant n-2c(\widetilde{G})+1$.

\textbf{Subcase 1.2.} $\widetilde{G}$ is a cycle-disjoint $U(\mathbb{Q})$-gain graph. Combining (1) and Subcase 1.1, we get $r(\widetilde{G})\geq n-2c(\widetilde{G})$.

\textbf{Case 2.} $p(\widetilde{G})\geq1$.

\textbf{Subcase 2.1.} $p(\widetilde{G})=1$.
Let $x$ be a pendant vertex of $\widetilde{G}$ and let $y$ be adjacent to $x$.
By Lemma \ref{le:2.5.}, we have $r(\widetilde{G})=r(\widetilde{G}-x-y)+2$.

If $d_{\widetilde{G}}(y)=2$, let $z$ be the other neighbour of $y$.
When $d_{\widetilde{G}}(z)=2$, we get $p(\widetilde{G}-x-y)=p(\widetilde{G})=1$ and $c(\widetilde{G}-x-y)=c(\widetilde{G})$. By applying induction hypothesis to $\widetilde{G}-x-y$, we obtain $r(\widetilde{G}-x-y)\geqslant n-2-2c(\widetilde{G}-x-y)-p(\widetilde{G}-x-y)+1$. Then $r(\widetilde{G})\geqslant n-2c(\widetilde{G})-p(\widetilde{G})+1$.
When $d_{\widetilde{G}}(z)\geq3$, we get $p(\widetilde{G}-x-y)=p(\widetilde{G})-1=0$ and $c(\widetilde{G}-x-y)=c(\widetilde{G})$. By applying induction hypothesis to $\widetilde{G}-x-y$, we have $r(\widetilde{G}-x-y)\geqslant n-2-2c(\widetilde{G}-x-y)$. Then $r(\widetilde{G})\geq n-2c(\widetilde{G})-p(\widetilde{G})+1$.

If $d_{\widetilde{G}}(y)\geqslant3$, let $\widetilde{H}_{1}, \widetilde{H}_{2}, \ldots, \widetilde{H}_{s}$ be the connected components of $\widetilde{G}-y$.
Without loss of generality, let $\widetilde{H}_{i}$ be the components containing 2-degree neighbors of $y$ for $i=1,...,r$,
$\widetilde{H}_{j}$ be the components containing  no 2-degree neighbors of $y$ for $j=r+1,...,s-1$, and let $H_{s}={x}$.
This arrangement implies that $\widetilde{H}_{i}$ has pendant vertices for $i=1,...,r$ and $\widetilde{H}_{j}$ has no pendant vertices for $j=r+1,...,s-1$. By the induction hypothesis, we get
$$r(\widetilde{H}_{i})\geq |V(\widetilde{H}_{i})|-2c(\widetilde{H}_{i})-p(\widetilde{H}_{i})+1~for~i=1,...,r;$$
$$r(\widetilde{H}_{j})\geq |V(\widetilde{H}_{j})|-2c(\widetilde{H}_{j})~for~j=r+1,...,s-1.$$
By Lemmas \ref{le:2.5.} and \ref{le:2.11.},
\begin{align}
r(\widetilde{G})&=r(\widetilde{G}-x-y)+2\nonumber
\\&=\sum^{r}_{i=1}r(\widetilde{H}_{i})+\sum^{s-1}_{j=r+1}r(\widetilde{H}_{j})+2\nonumber
\\&\geqslant \sum^{r}_{i=1}(|V(\widetilde{H}_{i})|-2c(\widetilde{H}_{i})-p(\widetilde{H}_{i})+1)+\sum^{s-1}_{j=r+1}(|V(\widetilde{H}_{j})|-2c(\widetilde{H}_{j}))+2\nonumber
\\&=n+r-2\sum^{s-1}_{i=1}c(\widetilde{H}_{i})-\sum^{r}_{i=1}p(\widetilde{H}_{i})\nonumber
\\&=n+r-2c(\widetilde{G}-x-y)-p(\widetilde{G}-x-y).\nonumber
\end{align}
Since $x$ is the unique pendant vertex of $\widetilde{G}$, $p(\widetilde{G}-x-y)=m$, where $m$ denotes the number of 2-degree neighbors of $y$. Since $c(\widetilde{H}_{s})=0$ and applying Lemma \ref{le:3.1.}$(3)$, we obtain
$$c(\widetilde{G}-x-y)=c(\widetilde{G}-y)=c(\widetilde{G})-d_{\widetilde{G}}(y)+s.$$
Then
\begin{equation}
r(\widetilde{G})\geq n-2c(\widetilde{G})+[2d_{\widetilde{G}}(y)+r-m-2s].
\end{equation}
By the fact that $d_{\widetilde{G}}(y)\geq s$ and Lemma \ref{le:3.1.}$(1)$, we obtain
\begin{equation}
2d_{\widetilde{G}}(y)+r-m-2s\geq0.
\end{equation}
Combining $(3)$ and $(4)$, we get $r(\widetilde{G})\geq n-2c(\widetilde{G})$, i.e., $r(\widetilde{G})\geq n-2c(\widetilde{G})-p(\widetilde{G})+1$ (since $p(\widetilde{G})=1$).

\textbf{Subcase 2.2.} $p(\widetilde{G})\geq2$.
Let $x$ be a pendant vertex and let $y$ be adjacent to $x$. As $p(\widetilde{G}-x)\geq1$, by applying induction hypothesis to $\widetilde{G}-x$, we obtain $r(\widetilde{G}-x)\geq n-1-2c(\widetilde{G}-x)-p(\widetilde{G}-x)+1.$ Obviously, $c(\widetilde{G}-x)=c(\widetilde{G})$.

If $d_{\widetilde{G}}(y)=2$, then $c(\widetilde{G}-x-y)=c(\widetilde{G})$ and
$1\leq p(\widetilde{G})-1\leq p(\widetilde{G}-x-y)\leq p(\widetilde{G})$. By using induction hypothesis to $\widetilde{G}-x-y$ and Lemma \ref{le:2.5.}, we have
\begin{align}
r(\widetilde{G})&=r(\widetilde{G}-x-y)+2\nonumber
\\&\geqslant n-2-2c(\widetilde{G}-x-y)-p(\widetilde{G}-x-y)+1+2\nonumber
\\&\geqslant n-2c(\widetilde{G})-p(\widetilde{G})+1.\nonumber
\end{align}

If $d_{\widetilde{G}}(y)\geq3$, then $p(\widetilde{G}-x)=p(\widetilde{G})-1$. By Lemma \ref{le:2.4.}, $r(\widetilde{G})\geq r(\widetilde{G}-x)\geq n-2c(\widetilde{G})-p(\widetilde{G})+1$.\quad $\square$

\section{ $U(\mathbb{Q})$-gain graphs $\widetilde{G}$ with $r(\widetilde{G})=|V(\widetilde{G})|-2c(\widetilde{G})$}

In this section,  we give a characterization of all $U(\mathbb{Q})$-gain graphs with $p(\widetilde{G})=0$ and $r(\widetilde{G})=|V(\widetilde{G})|-2c(\widetilde{G})$.

\noindent\begin{theorem}\label{th:5.1.}
Let $\widetilde{G}$ be a connected $U(\mathbb{Q})$-gain graph of order $n\geq2$ and $p(\widetilde{G})=0$. Then $r(\widetilde{G})=n-2c(\widetilde{G})$ if and only if $\widetilde{G}$ is a $U(\mathbb{Q})$-gain cycle which is of Type $1$.
\end{theorem}
\noindent\textbf{Proof.}
\textbf{Sufficiency:} Since $\widetilde{G}$ is a $U(\mathbb{Q})$-gain cycle $\widetilde{C}_{n}$ which is of Type $1$, by Lemma \ref{le:2.3.}, $r(\widetilde{G})=r(\widetilde{C}_{n})=n-2=n-2c(\widetilde{G})$.

\textbf{Necessity:} Since $p(\widetilde{G})=0$ and $r(\widetilde{G})=n-2c(\widetilde{G})$, by Theorem \ref{th:3.2.}, distinct cycles of $\widetilde{G}$ have no common vertices. Next, we will prove that $c(\widetilde{G})=1$.

If $c(\widetilde{G})\geq2$, then $\widetilde{G}$ has a pendant cycle $\widetilde{C}_{k}$. Let $x\in V(\widetilde{C}_{k})$ be a major vertex, $\widetilde{G}_{1}=\widetilde{G}-\widetilde{C}_{k}+x$ and let $\widetilde{G}_{2}=\widetilde{G}-\widetilde{C}_{k}$. Then $c(\widetilde{G}_{1})=c(\widetilde{G}_{2})=c(\widetilde{G})-1$, $p(\widetilde{G}_{1})=1$ and $p(\widetilde{G}_{2})\leq1$. By Theorem \ref{th:3.2.},
$$r(\widetilde{G}_{1})\geq |V(\widetilde{G}_{1})|-2c(\widetilde{G}_{1})-p(\widetilde{G}_{1})+1=n-k+3-2c(\widetilde{G}),$$
$$r(\widetilde{G}_{2})\geq |V(\widetilde{G}_{2})|-2c(\widetilde{G}_{2})=n-k+2-2c(\widetilde{G}).$$

If $\widetilde{C}_{k}$ is of Type $1$, then, by Lemma \ref{le:2.10.},
$r(\widetilde{G})=k-2+r(\widetilde{G}_{1})\geq n-2c(\widetilde{G})+1.$

If $\widetilde{C}_{k}$ is of Type $2$, then, by Lemma \ref{le:2.10.},
$r(\widetilde{G})=k+r(\widetilde{G}_{2})\geq n-2c(\widetilde{G})+2.$

If $\widetilde{C}_{k}$ is of Type $3$, then, by Lemma \ref{le:2.10.},
$r(\widetilde{G})\geq k-1+r(\widetilde{G}_{2})\geq n-2c(\widetilde{G})+1.$

If $\widetilde{C}_{k}$ is of Type $4$, then, by Lemma \ref{le:2.10.},
$r(\widetilde{G})=k-1+r(\widetilde{G}_{1})\geq n-2c(\widetilde{G})+2.$

So we get a contradiction. Thus $c(\widetilde{G})=1$ and $\widetilde{G}$ is a $U(\mathbb{Q})$-gain cycle which is of Type $1$ (by Lemma \ref{le:2.3.}).
\quad $\square$~

\begin{figure}[htbp]
  \centering
  \includegraphics[scale=1]{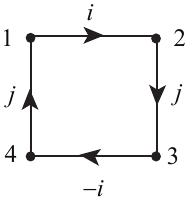}
  \caption{$\widetilde{G}$.}
\end{figure}

\noindent\begin{example}
Consider the $U(\mathbb{Q})$-gain graph $\widetilde{G}=\widetilde{C}_{4}$ in Fig. 1. Its adjacency matrix is given as follows:
\begin{align*}
 A(\widetilde{G})&=  \left (
 \begin{array}{ccccccccccccc}
 0 & i & 0 & -j\\
 -i & 0 & j & 0\\
 0 & -j & 0 & -i\\
 j & 0 & i & 0\\
 \end{array}
 \right).
 \end{align*}
We multiply $k$ on the left side of 1-th row and add it to 3-th row, and then we multiply $k$ on the left side of 2-th row and add it to 4-th row. Then:
\begin{align*}
r(\widetilde{G})=r\left(
\begin{array}{ccccccccccccc}
  0 & i & 0 & -j\\
 -i & 0 & j & 0\\
 0 & 0 & 0 & 0\\
 0 & 0 & 0 & 0\\
 \end{array}
 \right)
 =r \left(
\begin{array}{ccccccccccccc}
  0 & i & 0 & -j\\
 -i & 0 & j & 0\\
 \end{array}
 \right)
 =2.
 \end{align*}
Then $r(\widetilde{G})=2=|V(\widetilde{G})|-2c(\widetilde{G})$. Thus $\widetilde{G}$ is an extremal graph which satisfies the condition in Theorem \ref{th:5.1.}.
\end{example}

\section{ $U(\mathbb{Q})$-gain graphs $\widetilde{G}$ with $r(\widetilde{G})=|V(\widetilde{G})|-2c(\widetilde{G})+1$}

In this section, we will characterize all connected $U(\mathbb{Q})$-gain graphs with $p(\widetilde{G})=0$ and $r(\widetilde{G})=|V(\widetilde{G})|-2c(\widetilde{G})+1$.
There are some auxiliary results for the characterization of extremal $U(\mathbb{Q})$-gain graphs.

\noindent\begin{lemma}\label{le:4.7.}
Let $\widetilde{G}$ be a connected leaf-free $U(\mathbb{Q})$-gain graph with $c(\widetilde{G})\geq3$ and $r(\widetilde{G})=|V(\widetilde{G})|-2c(\widetilde{G})+1$. If $\widetilde{G}$ contains a pendant cycle $\widetilde{C}_{t}$, then $\widetilde{C}_{t}$ is of Type $1$.
\end{lemma}
\noindent\textbf{Proof.}
Let $x\in V(\widetilde{C}_{t})$ be a major vertex, $\widetilde{G}_{1}=\widetilde{G}-\widetilde{C}_{t}+x$ and let $\widetilde{G}_{2}=\widetilde{G}-\widetilde{C}_{t}$. Then $c(\widetilde{G}_{1})=c(\widetilde{G}_{2})=c(\widetilde{G})-1$, $p(\widetilde{G}_{1})=1$ and $p(\widetilde{G}_{2})\leq1$, by Theorem \ref{th:3.2.}, we obtain
$r(\widetilde{G}_{1})\geq |V(\widetilde{G}_{1})|-2c(\widetilde{G}_{1})-p(\widetilde{G}_{1})+1,$
$r(\widetilde{G}_{2})\geq |V(\widetilde{G}_{2})|-2c(\widetilde{G}_{2}).$

If $\widetilde{C}_{t}$ is of Type $2$, then, by Lemma \ref{le:2.10.}, $r(\widetilde{G})=t+r(\widetilde{G}_{2}).$
So $r(\widetilde{G}_{2})=|V(\widetilde{G}_{2})|-2c(\widetilde{G}_{2})-1$, a contradiction.

If $\widetilde{C}_{t}$ is of Type $3$, then, by Lemma \ref{le:2.10.}, $r(\widetilde{G})\geq t-1+r(\widetilde{G}_{2}).$
So $r(\widetilde{G}_{2})\leq |V(\widetilde{G}_{2})|-2c(\widetilde{G}_{2})$, thus $r(\widetilde{G}_{2})=|V(\widetilde{G}_{2})|-2c(\widetilde{G}_{2})$. By Theorem \ref{th:5.1.}, $\widetilde{G}_{2}$ is a $U(\mathbb{Q})$-gain cycle, a contradiction.

If $\widetilde{C}_{t}$ is of Type $4$, then, by Lemma \ref{le:2.10.}, $r(\widetilde{G})=t-1+r(\widetilde{G}_{1}).$
So $r(\widetilde{G}_{1})=|V(\widetilde{G}_{1})|-2c(\widetilde{G}_{1})-1$, a contradiction.

Thus $\widetilde{C}_{t}$ is of Type $1$.
\quad $\square$~\\

\noindent\begin{lemma}\label{le:4.3.}
Let $\widetilde{G}$ be a connected $U(\mathbb{Q})$-gain graph which is obtained from two $U(\mathbb{Q})$-gain graphs $\widetilde{H}$ and $\widetilde{K}$ by identifying the unique common vertex $v$, and $d_{\widetilde{K}}(v)\geq2$. If $r(\widetilde{K})\geq |V(\widetilde{K})|-2c(\widetilde{K})-p(\widetilde{K})+2$, then $r(\widetilde{G})\geq |V(\widetilde{G})|-2c(\widetilde{G})-p(\widetilde{G})+2$.
\end{lemma}
\noindent\textbf{Proof.}
\textbf{Case 1.} $d_{\widetilde{H}}(v)=1$.

\textbf{Subcase 1.1.} $\widetilde{H}$ is a $U(\mathbb{Q})$-gain path $\widetilde{P}_{m}$.
Then $c(\widetilde{G})=c(\widetilde{K})$ and $p(\widetilde{G})=p(\widetilde{K})+1$.
If $m$ is even, then, by Lemmas \ref{le:2.4.} and \ref{le:2.5.},
\begin{align}
r(\widetilde{G})&\geq r(\widetilde{K})+m-2\nonumber
\\&\geq |V(\widetilde{K})|-2c(\widetilde{K})-p(\widetilde{K})+2+m-2\nonumber
\\&=|V(\widetilde{G})|-2c(\widetilde{G})-p(\widetilde{G})+2.\nonumber
\end{align}
If $m$ is odd, then, by Lemma \ref{le:2.5.},
\begin{align}
r(\widetilde{G})&=r(\widetilde{K})+m-1\nonumber
\\&\geq |V(\widetilde{K})|-2c(\widetilde{K})-p(\widetilde{K})+2+m-1\nonumber
\\&=|V(\widetilde{G})|-2c(\widetilde{G})-p(\widetilde{G})+3.\nonumber
\end{align}

\textbf{Subcase 1.2.} $\widetilde{H}$ is not a $U(\mathbb{Q})$-gain path.
Let $u$ be a major vertex of $\widetilde{H}$ such that $d_{\widetilde{H}}(u,v)\leq d_{\widetilde{H}}(w,v)$ for any major vertex $w$ of $\widetilde{H}$ and let $\widetilde{P}_{m}$ be the path with $v, u$ as its end vertices. Let $\widetilde{H}_{1}=\widetilde{H}-\widetilde{P}_{m}+u$.

If $\widetilde{H}_{1}$ is a $U(\mathbb{Q})$-gain cycle, let $\widetilde{G}_{1}=\widetilde{K}+(\widetilde{P}_{m}-v)$ and let $\widetilde{G}_{2}=\widetilde{K}+(\widetilde{P}_{m}-v-u)$. Then $c(\widetilde{G})=c(\widetilde{K})+1$ and $p(\widetilde{G})=p(\widetilde{K})$.
When $m$ is even, by Lemmas \ref{le:2.4.} and \ref{le:2.5.}, we get
$r(\widetilde{G}_{1})\geq r(\widetilde{K})+m-2$ and $r(\widetilde{G}_{2})=r(\widetilde{K})+m-2$.
By Lemma \ref{le:2.10.} and Theorem \ref{th:3.2.}, we obtain
\begin{align}
r(\widetilde{G})&\geq |V(\widetilde{H}_{1})|-2+r(\widetilde{G}_{1})\nonumber
\\&\geq |V(\widetilde{H}_{1})|-2+(r(\widetilde{K})+m-2)\nonumber
\\&\geq |V(\widetilde{H}_{1})|-2+(|V(\widetilde{K})|-2c(\widetilde{K})-p(\widetilde{K})+2)+m-2\nonumber
\\&=|V(\widetilde{G})|-2c(\widetilde{G})-p(\widetilde{G})+2.\nonumber
\end{align}
When $m$ is odd, by Lemmas \ref{le:2.4.} and \ref{le:2.5.},
$r(\widetilde{G}_{1})=r(\widetilde{K})+m-1$ and $r(\widetilde{G}_{2})\geq r(\widetilde{K})+m-3$.
By Lemma \ref{le:2.10.} and Theorem \ref{th:3.2.},
\begin{align}
r(\widetilde{G})&\geq |V(\widetilde{H}_{1})|-1+r(\widetilde{G}_{2})\nonumber
\\&\geq |V(\widetilde{H}_{1})|-1+(r(\widetilde{K})+m-3)\nonumber
\\&\geq |V(\widetilde{H}_{1})|-1+(|V(\widetilde{K})|-2c(\widetilde{K})-p(\widetilde{K})+2)+m-3\nonumber
\\&=|V(\widetilde{G})|-2c(\widetilde{G})-p(\widetilde{G})+2.\nonumber
\end{align}

If $\widetilde{H}_{1}$ is not a $U(\mathbb{Q})$-gain cycle, then $r(\widetilde{H}_{1})\geq |V(\widetilde{H}_{1})|-2c(\widetilde{H}_{1})-p(\widetilde{H}_{1})+1$ (by Theorems \ref{th:3.2.} and \ref{th:5.1.}), $c(\widetilde{G})=c(\widetilde{K})+c(\widetilde{H}_{1})$ and $p(\widetilde{G})=p(\widetilde{K})+p(\widetilde{H}_{1})$.
By Lemma \ref{le:4.2.}, we get
\begin{align}
r(\widetilde{G})&\geq r(\widetilde{K})+r(\widetilde{H}_{1})+m-3\nonumber
\\&\geq (|V(\widetilde{K})|-2c(\widetilde{K})-p(\widetilde{K})+2)+(|V(\widetilde{H}_{1})|-2c(\widetilde{H}_{1})-p(\widetilde{H}_{1})+1)+m-3\nonumber
\\&=|V(\widetilde{G})|-2c(\widetilde{G})-p(\widetilde{G})+2.\nonumber
\end{align}

\textbf{Case 2.} $d_{\widetilde{H}}(v)\geq2$.

If $\widetilde{H}-v$ is a $U(\mathbb{Q})$-gain cycle, then $c(\widetilde{K})\leq c(\widetilde{G})-2$, $p(\widetilde{K})=p(\widetilde{G})$ and $r(\widetilde{H}-v)\geq |V(\widetilde{H}-v)|-2$. By Lemma \ref{le:4.1.}, we obtain
\begin{align}
r(\widetilde{G})&=r(\widetilde{K})+r(\widetilde{H}-v)-1\nonumber
\\&\geq (|V(\widetilde{K})|-2c(\widetilde{K})-p(\widetilde{K})+2)+(|V(\widetilde{H}-v)|-2)-1\nonumber
\\&=|V(\widetilde{G})|-2c(\widetilde{G})-p(\widetilde{G})+3.\nonumber
\end{align}
If $\widetilde{H}-v$ is not a $U(\mathbb{Q})$-gain cycle, then $c(\widetilde{G})=c(\widetilde{K})+c(\widetilde{H}-v)+d_{\widetilde{H}}(v)-1$ (by Lemma \ref{le:3.1.}), $p(\widetilde{H}-v)+p(\widetilde{K})\leq p(\widetilde{G})+d_{\widetilde{H}}(v)$. By Theorem \ref{th:3.2.}, $r(\widetilde{H}-v)\geq |V(\widetilde{H}-v)|-2c(\widetilde{H}-v)-p(\widetilde{H}-v)+1$. By Lemma \ref{le:4.1.}, we obtain
\begin{align}
r(\widetilde{G})&=r(\widetilde{K})+r(\widetilde{H}-v)-1\nonumber
\\&\geq (|V(\widetilde{K})|-2c(\widetilde{K})-p(\widetilde{K})+2)+(|V(\widetilde{H}-v)|-2c(\widetilde{H}-v)-p(\widetilde{H}-v)+1)-1\nonumber
\\&\geq |V(\widetilde{G})|-2c(\widetilde{G})-p(\widetilde{G})+2.~~\square\nonumber
\end{align}

If a vertex $v$ is covered by all maximum matchings in a $U(\mathbb{Q})$-gain graph, then $v$ is called a \emph{covered vertex}.
By Lemma $3.1$ in \cite{ccz} and Lemma \ref{le:2.0.}, we have the following theorem for $U(\mathbb{Q})$-gain trees.
An \emph{internal path} of $\widetilde{G}$ is a path with every inner vertex (except end vertices) of degree $2$ in $\widetilde{G}$.
\noindent\begin{theorem}\label{th:5.4.}
Let $\widetilde{T}$ be a $U(\mathbb{Q})$-gain tree of order $n$ with $p(\widetilde{T})\geq3$. Then $r(\widetilde{T})=n-p(\widetilde{T})+1$ if and only if $\widetilde{T}$ satisfies the following two conditions:
\begin{enumerate}[(1)]
  \item All internal paths from a leaf to a major vertex are odd;
  \item Let $\widetilde{P}_{k}$ be an internal path from a leaf $u$ of $\widetilde{T}$ to a major vertex $v$. Then $r(\widetilde{T}_{1})=|V(\widetilde{T}_{1})|-p(\widetilde{T}_{1})+1$, where $\widetilde{T}_{1}=\widetilde{T}-\widetilde{P}_{k}+v$ and $v$ is a covered vertex of $\widetilde{T}_{1}$.
\end{enumerate}
\end{theorem}

\noindent\begin{lemma}\label{le:4.4.}
Let $\widetilde{G}$ be a connected $U(\mathbb{Q})$-gain graph with a unique $U(\mathbb{Q})$-gain cycle $\widetilde{C}_{m}$ and $p(\widetilde{G})\geq1$. Then $r(\widetilde{G})=|V(\widetilde{G})|-2c(\widetilde{G})-p(\widetilde{G})+1$ if and only if $\widetilde{C}_{m}$ is of Type $1$ and $\widetilde{G}$ is obtained from a $U(\mathbb{Q})$-gain tree $\widetilde{T}$ with $r(\widetilde{T})=|V(\widetilde{T})|-p(\widetilde{T})+1$, by attaching $\widetilde{C}_{m}$ at a leaf $x$ of $\widetilde{T}$.
\end{lemma}
\noindent\textbf{Proof.}
\textbf{Sufficiency:} By Lemma \ref{le:2.10.},
$$r(\widetilde{G})=m-2+r(\widetilde{T})=m-2+|V(\widetilde{T})|-p(\widetilde{T})+1=|V(\widetilde{G})|-2c(\widetilde{G})-p(\widetilde{G})+1.$$

\textbf{Necessity:} Since $\widetilde{G}$ has a unique $U(\mathbb{Q})$-gain cycle $\widetilde{C}_{m}$, then $\widetilde{C}_{m}$ is a block, each major vertex of $\widetilde{G}$ on $\widetilde{C}_{m}$ must be a cut-point of $\widetilde{G}$. Let $x$ be a cut-point of $\widetilde{G}$ which lies on $\widetilde{C}_{m}$. Let $\widetilde{H}_{1}$ be a component of $\widetilde{G}-x$ such that $V(\widetilde{H}_{1})\cap V(\widetilde{C}_{m})=\emptyset$, $\widetilde{H}=\widetilde{H}_{1}+x$ and let $\widetilde{G}_{1}=\widetilde{C}_{m}+\widetilde{H}_{1}$. Next, we will prove that $\widetilde{C}_{m}$ is of Type $1$.

If $\widetilde{C}_{m}$ is of Type $2$ or $3$, then, by Lemma \ref{le:2.3.}, $r(\widetilde{C}_{m})=m=m-2c(\widetilde{C}_{m})-p(\widetilde{C}_{m})+2$.
Thus by Lemma \ref{le:4.3.}, $r(\widetilde{G})=|V(\widetilde{G})|-2c(\widetilde{G})-p(\widetilde{G})+2$, a contradiction.
If $\widetilde{C}_{m}$ is of Type $4$. Since $c(\widetilde{H})=c(\widetilde{G}_{1})-1$ and $p(\widetilde{H})=p(\widetilde{G}_{1})+1$, by Lemma \ref{le:2.10.} and Theorem \ref{th:3.2.},
\begin{align}
r(\widetilde{G}_{1})&=m-1+r(\widetilde{H})\nonumber
\\&\geq m-1+(|V(\widetilde{H})|-2c(\widetilde{H})-p(\widetilde{H})+1)\nonumber
\\&=|V(\widetilde{G}_{1})|-2c(\widetilde{G}_{1})-p(\widetilde{G}_{1})+2.\nonumber
\end{align}
Thus by Lemma \ref{le:4.3.}, $r(\widetilde{G})=|V(\widetilde{G})|-2c(\widetilde{G})-p(\widetilde{G})+2$, a contradiction. Hence, $\widetilde{C}_{m}$ is of Type $1$.
For the rest of the proof, we consider two cases.

\textbf{Case 1.} $\widetilde{C}_{m}$ has a unique major vertex $x$ of $\widetilde{G}$.

In this case,  $\widetilde{G}-\widetilde{C}_{m}+x$ is a $U(\mathbb{Q})$-gain tree $\widetilde{T}$.

If $d_{\widetilde{G}}(x)=3$, then $\widetilde{C}_{m}$ is a pendant cycle of $\widetilde{G}$. By Lemma \ref{le:2.10.}, $r(\widetilde{G})=m-2+r(\widetilde{T})$, then $r(\widetilde{T})=|V(\widetilde{T})|-p(\widetilde{T})+1$. Thus $\widetilde{G}$ is obtained from $\widetilde{T}$ with $r(\widetilde{T})=|V(\widetilde{T})|-p(\widetilde{T})+1$, by attaching $\widetilde{C}_{m}$ at a leaf $x$ of $\widetilde{T}$.

If $d_{\widetilde{G}}(x)\geq4$, by Lemma \ref{le:2.10.}, $r(\widetilde{G})=m-2+r(\widetilde{T})$, then $r(\widetilde{T})=|V(\widetilde{T})|-p(\widetilde{T}_{1})$, a contradiction.

\textbf{Case 2.} $\widetilde{C}_{m}$ has at least two major vertices of $\widetilde{G}$.

In this case, we will prove that $r(\widetilde{G})\geq|V(\widetilde{G})|-2c(\widetilde{G})-p(\widetilde{G})+2$.
By Lemma \ref{le:4.3.}, we may assume that $\widetilde{C}_{m}$ has exactly two major vertices of $\widetilde{G}$, say $x$ and $y$.
Let $\widetilde{T}_{x}$, $\widetilde{T}_{y}$ be the pendant trees attached at $x$, $y$ of $\widetilde{C}_{m}$ and let $\widetilde{G}_{1}=\widetilde{G}-(\widetilde{T}_{x}-x)$, $\widetilde{G}_{2}=\widetilde{G}-(\widetilde{T}_{y}-y)$.
By Lemma \ref{le:4.3.}, $r(\widetilde{G}_{1})=|V(\widetilde{G}_{1})|-2c(\widetilde{G}_{1})-p(\widetilde{G}_{1})+1$.
Then $r(\widetilde{T}_{y})=|V(\widetilde{T}_{y})|-p(\widetilde{T}_{y})+1$ and $\widetilde{C}_{m}$ is a pendant cycle of $\widetilde{G}_{1}$ (by Case 1). In a similar way, $r(\widetilde{T}_{x})=|V(\widetilde{T}_{x})|-p(\widetilde{T}_{x})+1$ and $\widetilde{C}_{m}$ is a pendant cycle of $\widetilde{G}_{2}$. By Theorem \ref{th:5.4.}, $r(\widetilde{T}_{x}-x)\geq|V(\widetilde{T}_{x}-x)|-p(\widetilde{T}_{x}-x)+2$ and $r(\widetilde{G}_{1}-x)\geq|V(\widetilde{G}_{1}-x)|-p(\widetilde{G}_{1}-x)+2$. By Lemma \ref{le:2.4.}, we obtain
$$r(\widetilde{G})\geq r(\widetilde{G}-x)
=r(\widetilde{T}_{x}-x)+r(\widetilde{G}_{1}-x)
\geq|V(\widetilde{G})|-2c(\widetilde{G})-p(\widetilde{G})+2,$$
a contradiction.
\quad $\square$~\\

For a simple graph, let $C_{p}$ and $C_{q}$ be two vertex-disjoint cycles with $v\in V(C_{p})$ and $u\in V(C_{q})$,
and let $P_{l}=v_{1}v_{2}\cdots v_{l}~(l\geq1)$ be a  path of length $l-1$.
Let $\infty(p,l,q)$ (as shown in Fig. 2) be the graph obtained from $C_{p}$, $C_{q}$ and $P_{l}$ by identifying $v$ with $v_{1}$ and $u$ with $v_{l}$, respectively.  When $l=1$, the graph $\infty(p,1,q)$ (as shown in Fig. 2) is obtained from $C_{p}$ and $C_{q}$ by identifying $v$ with $u$.

Let $P_{p+2},P_{l+2},P_{q+2}$ be three paths, where min$\{p,l,q\}\geq0$ and at most one of $p,l,q$ is 0. Let $\theta(p,l,q)$ (as shown in Fig. 2) be the graph obtained from $P_{p+2}$, $P_{l+2}$ and $P_{q+2}$ by identifying the three initial vertices and terminal vertices.

For a $U(\mathbb{Q})$-gain graph, let $\widetilde{\infty}(p,l,q)=(\infty(p,l,q),U(\mathbb{Q}),\varphi')$ and $\widetilde{\theta}(p,l,q)=(\theta(p,l,q),U(\mathbb{Q}),\\\varphi'')$.

\begin{figure}[htbp]
 \centering
 \includegraphics[scale=0.75]{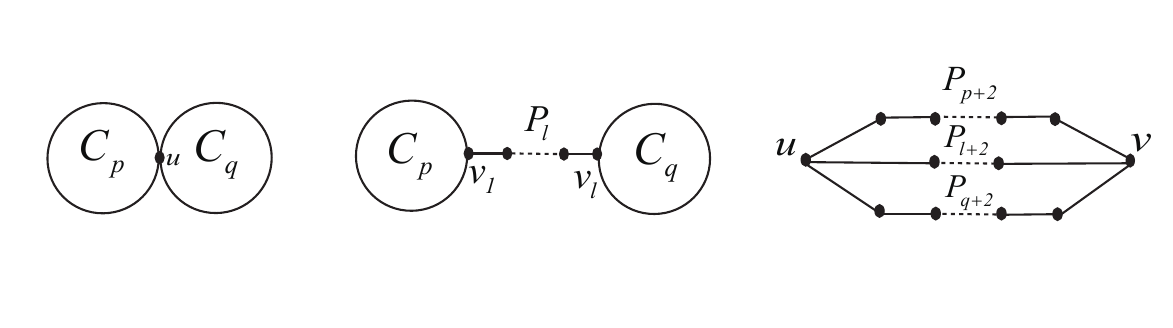}
 \caption{$\infty(p,1,q)$, $\infty(p,l,q)$ and $\theta(p,l,q)$.}
\end{figure}

\noindent\begin{lemma}\label{le:4.5.}
Let $\widetilde{G}$ be a connected $U(\mathbb{Q})$-gain graph of order $n$ with $r(\widetilde{G})=n-2c(\widetilde{G})-p(\widetilde{G})+1$ and let $\widetilde{C}_{m}$ be a block in $\widetilde{G}$. Then
\begin{enumerate}[(1)]
  \item $\widetilde{C}_{m}$ is of Type $1$;
  \item Either $\widetilde{C}_{m}$ is a pendant cycle of $\widetilde{G}$, or $\widetilde{G}=\widetilde{\infty}(m,1,n-m+1)$ such that each cycle of $\widetilde{G}$ is of Type $1$.
\end{enumerate}
\end{lemma}
\noindent\textbf{Proof.}
Since $\widetilde{C}_{m}$ is a block, each major vertex of $\widetilde{G}$ on $\widetilde{C}_{m}$ must be a cut-point of $\widetilde{G}$.
Let $x$ be a cut-point of $\widetilde{G}$ which lies on $\widetilde{C}_{m}$.
Let $\widetilde{H}_{1}$ be a component of $\widetilde{G}-x$ such that $V(\widetilde{H}_{1})\cap V(\widetilde{C}_{m})=\emptyset$ and  let $\widetilde{H}=\widetilde{H}_{1}+x$. Next, we will prove that $\widetilde{C}_{m}$ is of Type $1$.

If $\widetilde{C}_{m}$ is of Type $2$ or $3$, then, by Lemma \ref{le:2.3.}, $r(\widetilde{C}_{m})=m=m-2c(\widetilde{C}_{m})-p(\widetilde{C}_{m})+2$, thus $r(\widetilde{G})\geq n-2c(\widetilde{G})-p(\widetilde{G})+2$ (by Lemma \ref{le:4.3.}), a contradiction.

If $\widetilde{C}_{m}$ is of Type $4$. Let $\widetilde{G}_{1}=\widetilde{C}_{m}+\widetilde{H}_{1}$.
When $\widetilde{H}$ is a $U(\mathbb{Q})$-gain cycle which is of Type 1, by Lemma \ref{le:2.10.}, $r(\widetilde{G}_{1})=|V(\widetilde{H})|-2+r(\widetilde{C}_{m})=|V(\widetilde{G}_{1})|-2=|V(\widetilde{G}_{1})|-2c(\widetilde{G}_{1})-p(\widetilde{G}_{1})+2$, thus $r(\widetilde{G})\geq n-2c(\widetilde{G})-p(\widetilde{G})+2$ (by Lemma \ref{le:4.3.}), a contradiction.
When $\widetilde{H}$ is not a $U(\mathbb{Q})$-gain cycle or a cycle which is not of Type 1, then $c(\widetilde{H})=c(\widetilde{G}_{1})-1$ and $p(\widetilde{H})\leq p(\widetilde{G}_{1})+1$. By Theorems \ref{th:3.2.} and \ref{th:5.1.}, $r(\widetilde{H})\geq |V(\widetilde{H})|-2c(\widetilde{H})-p(\widetilde{H})+1$. By Lemma \ref{le:2.10.} and Theorem \ref{th:3.2.},
$$r(\widetilde{G}_{1})=m-1+r(\widetilde{H})\geq m-1+(|V(\widetilde{H})|-2c(\widetilde{H})-p(\widetilde{H})+1)\geq |V(\widetilde{G}_{1})|-2c(\widetilde{G}_{1})-p(\widetilde{G}_{1})+2.$$
Thus $r(\widetilde{G})\geq n-2c(\widetilde{G})-p(\widetilde{G})+2$ (by Lemma \ref{le:4.3.}), a contradiction. Thus (1) holds.

The proof of (2) is given by two cases as follows.

\textbf{Case 1.} $\widetilde{C}_{m}$ has a unique major vertex $x$ of $\widetilde{G}$.

If $d_{\widetilde{G}}(x)=3$, then $\widetilde{C}_{m}$ is a pendant cycle of $\widetilde{G}$.
If $d_{\widetilde{G}}(x)\geq4$, we assume $\widetilde{G}_{1}=\widetilde{G}-\widetilde{C}_{m}+x$.
By Lemma \ref{le:2.10.}, $r(\widetilde{G})=m-2+r(\widetilde{G}_{1})$.
Since $r(\widetilde{G})=n-2c(\widetilde{G})-p(\widetilde{G})+1$, $r(\widetilde{G}_{1})=|V(\widetilde{G}_{1})|-2c(\widetilde{G}_{1})-p(\widetilde{G}_{1})$, by Theorems \ref{th:3.2.} and \ref{th:5.1.}, $\widetilde{G}_{1}$ is a $U(\mathbb{Q})$-gain cycle which is of Type $1$.
Thus $\widetilde{G}=\widetilde{\infty}(m,1,n-m+1)$ and each cycle of $\widetilde{G}$ is of Type $1$.

\textbf{Case 2.} $\widetilde{C}_{m}$ has at least two major vertices of $\widetilde{G}$.

By Lemma \ref{le:4.3.}, we just need to prove that $\widetilde{C}_{m}$ has two major vertices of $\widetilde{G}$. We will prove that if $\widetilde{C}_{m}$ has two major vertices of $\widetilde{G}$, then $r(\widetilde{G})\geq n-2c(\widetilde{G})-p(\widetilde{G})+2$. We proceed by induction on $c(\widetilde{G})$. If $c(\widetilde{G})=1$, the contradiction is deduced by Lemma \ref{le:4.4.}. Now we assume that the conclusion holds for any connected $U(\mathbb{Q})$-gain graph with $c-1(\geq1)$ cycles, while $c(\widetilde{G})=c$. Let $\widetilde{C}$ be a cycle of $\widetilde{G}$ distinct to $\widetilde{C}_{m}$.
If $\widetilde{C}$ has only one major vertex $x$ of $\widetilde{G}$, then $\widetilde{C}$ is a pendant cycle of $\widetilde{G}$ and $\widetilde{C}$ is of Type 1 (by Case 1). Let $\widetilde{G}_{1}=\widetilde{G}-\widetilde{C}+x$, then $c(\widetilde{G})=c(\widetilde{G}_{1})+1$ and $p(\widetilde{G})=p(\widetilde{G}_{1})-1$. The induction hypothesis imply that  $r(\widetilde{G}_{1})\geq|V(\widetilde{G}_{1})|-2c(\widetilde{G}_{1})-p(\widetilde{G}_{1})+2$. By Lemma \ref{le:2.10.},
$$r(\widetilde{G})=|V(\widetilde{C})|-2+r(\widetilde{G}_{1})\geq n-2c(\widetilde{G})-p(\widetilde{G})+2,$$
as required.
If $\widetilde{C}$ has at least two major vertices.
Let $x, y$ be two major vertices of $\widetilde{C}$ and
let $\widetilde{P}_{m}$ be the path contained in $\widetilde{C}$ with $x$ and $y$ as its end vertices (internal vertices (if any) of $\widetilde{P}$ are not major vertices).
When $m=2$, let $\widetilde{G}_{1}=\widetilde{G}-e_{xy}$. By Lemma \ref{le:2.6.} and the induction hypothesis, we get
$$r(\widetilde{G})\geq r(\widetilde{G}-e_{xy})-2\geq|V(\widetilde{G}_{1})|-2c(\widetilde{G}_{1})-p(\widetilde{G}_{1})+2-2=n-2c-p(\widetilde{G})+2.$$
When $m=3$, let $V(\widetilde{P}_{m})=\{x,z,y\}$ and let $\widetilde{G}_{1}=\widetilde{G}-z$. By Lemma \ref{le:2.4.} and the induction hypothesis, we have
$$r(\widetilde{G})\geq r(\widetilde{G}-z)\geq|V(\widetilde{G}_{1})|-2c(\widetilde{G}_{1})-p(\widetilde{G}_{1})+2=n-2c-p(\widetilde{G})+3.$$
When $m\geq4$, let $z\in V(\widetilde{P}_{m})$ be adjacent to $x$ and let $\widetilde{G}_{1}=\widetilde{G}-z$. By Lemma \ref{le:2.4.} and the induction hypothesis, we get
$$r(\widetilde{G})\geq r(\widetilde{G}-z)\geq|V(\widetilde{G}_{1})|-2c(\widetilde{G}_{1})-p(\widetilde{G}_{1})+2=n-2c-p(\widetilde{G})+2,$$
showing that $r(\widetilde{G})\geq n-2c-p(\widetilde{G})+2$, a contradiction. Hence, $\widetilde{C}_{m}$ has a major vertex of $\widetilde{G}$ (i.e., Case 1).

Thus (2) holds.
\quad $\square$~\\

\noindent\begin{lemma}\label{le:4.6.}
Let $\widetilde{G}$ be a connected $U(\mathbb{Q})$-gain graph of order $n$ with $r(\widetilde{G})=n-2c(\widetilde{G})-p(\widetilde{G})+1$. Then each cycle (if any) of $\widetilde{G}$ is of Type $1$.
\end{lemma}
\noindent\textbf{Proof.}
This theorem is equivalent to proving that if $\widetilde{G}$ contains a $U(\mathbb{Q})$-gain cycle $\widetilde{C}_{m}$ which is not of Type 1, then $r(\widetilde{G})\geq n-2c(\widetilde{G})-p(\widetilde{G})+2$.

Now, we proceed by induction on $c(\widetilde{G})$.
If $c(\widetilde{G})=1$, the conclusion is proved by Lemma \ref{le:4.4.}.
Suppose the result holds for graphs with $c-1(\geq 1)$ cycles, while $c(\widetilde{G})=c$.
If $\widetilde{C}_{m}$ is a block of $\widetilde{G}$, the conclusion is proved by Lemma \ref{le:4.5.}.
Now assume that $\widetilde{C}_{m}$ is not a block of $\widetilde{G}$.
Let $\widetilde{C}$ be a $U(\mathbb{Q})$-gain cycle such that $\widetilde{C}$ and $\widetilde{C}_{m}$ share at least two common vertices
and let $x\in V(\widetilde{C}), x \notin V(\widetilde{C}_{m})$ be a vertex adjacent to a common vertex $y$ of $\widetilde{C}$ and $\widetilde{C}_{m}$.
If $d_{\widetilde{C}}(x)=2$, then $c(\widetilde{G}-x)\leq c(\widetilde{G})-1$ and $p(\widetilde{G}-x)\leq p(\widetilde{G})+1$. Lemma \ref{le:2.4.} and the induction hypothesis imply that
$$r(\widetilde{G})\geq r(\widetilde{G}-x)\geq|V(\widetilde{G}-x)|-2c(\widetilde{G}-x)-p(\widetilde{G}-x)+2\geq n-2c-p(\widetilde{G})+2,$$
as required.
If $d_{\widetilde{C}}(x)\geq3$, then $c(\widetilde{G}-e_{xy})\leq c(\widetilde{G})-1$ and $p(\widetilde{G}-e_{xy})=p(\widetilde{G})$. Lemma \ref{le:2.6.} and the induction hypothesis imply that
$$r(\widetilde{G})\geq r(\widetilde{G}-e_{xy})-2\geq|V(\widetilde{G}-e_{xy})|-2c(\widetilde{G}-e_{xy})-p(\widetilde{G}-e_{xy})+2-2\geq n-2c-p(\widetilde{G})+2,$$
as required.
\quad $\square$~\\

\noindent\begin{theorem}\label{th:5.2.}
Let $\widetilde{G}$ be a connected $U(\mathbb{Q})$-gain graph of order $n$ with $p(\widetilde{G})=0$ and $c(\widetilde{G})=2$. Then $r(\widetilde{G})=n-2c(\widetilde{G})+1=n-3$ if and only if $\widetilde{G}=\widetilde{\infty}(p,l,q)$, where each cycle of $\widetilde{G}$ is of Type $1$  and $l$ is odd, or $\widetilde{G}=\widetilde{\theta}(p',l',q') $, where each cycle of $\widetilde{G}$ is of Type $1$ and $p', l', q'$ are odd.
\end{theorem}
\noindent\textbf{Proof.}
\textbf{Sufficiency:}
Assume $\widetilde{G}=\widetilde{\infty}(p,l,q)$, where each cycle of $\widetilde{G}$ is of Type $1$ and $l$ is odd.
Let $v$ be the $3$-degree vertex on the $\widetilde{C}_{p}$ and let $\widetilde{G}_{1}=\widetilde{G}-\widetilde{C}_{p}+v$. By Lemma \ref{le:2.5.}, $r(\widetilde{G}_{1})=r(\widetilde{C}_{q})+l-1=q+l-3$. By Lemma \ref{le:2.10.}, $r(\widetilde{G})=p-2+r(\widetilde{G}_{1})=n-3$.

Let $\widetilde{G}=\widetilde{\theta}(p',l',q')$, where each cycle of $\widetilde{G}$ is of Type $1$ and $p',l',q'$ are odd.
When $p',l',q'\leq3$, by calculations, we can know $r(\widetilde{G})=n-3$.
When at least one of $p',l',q'$ is greater than or equal to 3, by Lemma \ref{le:2.7.} $k$ times, $r(\widetilde{G})=r(\widetilde{G}_{1})+4k$, where $\widetilde{G}_{1}=\widetilde{\theta}(3,3,3),~\widetilde{\theta}(3,3,1),~\widetilde{\theta}(3,1,1)~or~\widetilde{\theta}(1,1,1)$. Then  $r(\widetilde{G})=n-3$.

\textbf{Necessity:} Let $\widetilde{G}$ be a $U(\mathbb{Q})$-gain graph with $p(\widetilde{G})=0$, $c(\widetilde{G})=2$ and $r(\widetilde{G})=n-2c(\widetilde{G})-p(\widetilde{G})+1=n-3$. Then $\widetilde{G}=\widetilde{\infty}(p,l,q)$, where each cycle of $\widetilde{G}$ is of Type $1$, or $\widetilde{G}=\widetilde{\theta}(p',l',q')$, where each cycle of $\widetilde{G}$ is of Type $1$ (by Lemma \ref{le:4.6.}).

\textbf{Case 1.} $\widetilde{G}=\widetilde{\infty}(p,l,q)$, where each cycle of $\widetilde{G}$ is of Type $1$.

Let $u$ be the $3$-degree vertex on the $\widetilde{C}_{p}$, $v$ be the $3$-degree vertex on the $\widetilde{C}_{q}$ and let $\widetilde{G}_{1}=\widetilde{G}-\widetilde{C}_{p}+u$. If $l$ is even, by Lemmas \ref{le:2.1.} and \ref{le:2.5.}, $r(\widetilde{G}_{1})=r(\widetilde{C}_{q}-v)+l=l+q-2$. By Lemma \ref{le:2.10.}, $r(\widetilde{G})=p-2+r(\widetilde{G}_{1})=n-2$, a contradiction; thus $l$ is odd.

\textbf{Case 2.} $\widetilde{G}=\widetilde{\theta}(p',l',q')$, where each cycle of $\widetilde{G}$ is of Type $1$.

If $p'$ is even, then $l', q'$ are even. Let $u$ be a 3-degree vertex, by Lemmas \ref{le:2.4.} and \ref{le:2.5.}, $r(\widetilde{G})\geq r(\widetilde{G}-u)=n-2$, a contradiction; thus $p',l',q'$ are odd.
\quad $\square$~\\

\begin{figure}[htbp]
  \centering
  \includegraphics[scale=0.9]{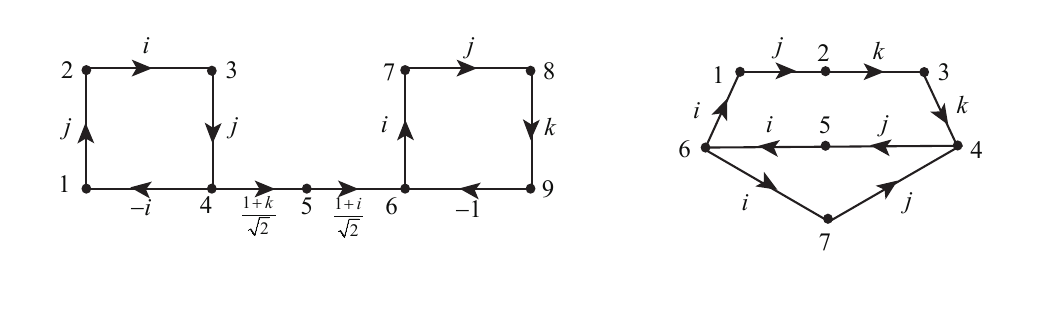}
  \caption{$\widetilde{G}_{1}$ and $\widetilde{G}_{2}$.}
\end{figure}

\noindent\begin{example}
Consider the $U(\mathbb{Q})$-gain graphs $\widetilde{G}_{1}$ and $\widetilde{G}_{2}$ in Fig. 3. Their adjacency matrices are given as follows:
\begin{align*}
 A(\widetilde{G}_{1})&= \left (
\begin{array}{ccccccccccccc}
 0 & j & 0 & i & 0 & 0 & 0 & 0 & 0\\
 -j & 0 & i & 0 & 0 & 0 & 0 & 0 & 0\\
 0 & -i & 0 & j & 0 & 0 & 0 & 0 & 0\\
 -i & 0 & -j & 0 & (1+k)/\sqrt{2} & 0 & 0 & 0 & 0\\
 0 & 0 & 0 & (1-k)/\sqrt{2} & 0 & (1+i)/\sqrt{2} & 0 & 0 & 0\\
 0 & 0 & 0 & 0 & (1-i)/\sqrt{2} & 0 & i & 0 & -1  \\
 0 & 0 & 0 & 0 & 0 & -i & 0 & j & 0\\
 0 & 0 & 0 & 0 & 0 & 0 & -j & 0 & k\\
 0 & 0 & 0 & 0 & 0 & -1 & 0 & -k & 0\\
\end{array}
 \right).
\end{align*}
We multiply $-k$ on the left side of 1-th row and add it to 3-th row,
multiply $k$ on the left side of 2-th row and add it to 4-th row,
multiply $-k$ on the left side of 8-th row and add it to 6-th row,
multiply $-i$ on the left side of 9-th row and add it to 7-th row:
\begin{align*}
  r(\widetilde{G}_{1})=r \left(
\begin{array}{ccccccccccccc}
 0 & j & 0 & i & 0 & 0 & 0 & 0 & 0\\
 -j & 0 & i & 0 & 0 & 0 & 0 & 0 & 0\\
 0 & 0 & 0 & 0 & 0 & 0 & 0 & 0 & 0\\
 0 & 0 & 0 & 0 & (1+k)/\sqrt{2} & 0 & 0 & 0 & 0\\
 0 & 0 & 0 & (1-k)/\sqrt{2} & 0 & (1+i)/\sqrt{2} & 0 & 0 & 0\\
 0 & 0 & 0 & 0 & (1-i)/\sqrt{2} & 0 & 0 & 0 & 0  \\
 0 & 0 & 0 & 0 & 0 & 0 & 0 & 0 & 0\\
 0 & 0 & 0 & 0 & 0 & 0 & -j & 0 & k\\
 0 & 0 & 0 & 0 & 0 & 1 & 0 & -k & 0\\
 \end{array}
 \right).
 \end{align*}
We multiply $-(1-i)/(1+k)$ on the left side of 4-th row and add it to 6-th row:
 \begin{align*}
  r(\widetilde{G}_{1})&=r \left(
\begin{array}{ccccccccccccc}
 0 & j & 0 & i & 0 & 0 & 0 & 0 & 0\\
 -j & 0 & i & 0 & 0 & 0 & 0 & 0 & 0\\
 0 & 0 & 0 & 0 & 0 & 0 & 0 & 0 & 0\\
 0 & 0 & 0 & 0 & (1+k)/\sqrt{2} & 0 & 0 & 0 & 0\\
 0 & 0 & 0 & (1-k)/\sqrt{2} & 0 & (1+i)/\sqrt{2} & 0 & 0 & 0\\
 0 & 0 & 0 & 0 & 0 & 0 & 0 & 0 & 0  \\
 0 & 0 & 0 & 0 & 0 & 0 & 0 & 0 & 0\\
 0 & 0 & 0 & 0 & 0 & 0 & -j & 0 & k\\
 0 & 0 & 0 & 0 & 0 & 1 & 0 & -k & 0\\
 \end{array}
 \right)\\
 &=r \left(
\begin{array}{ccccccccccccc}
 0 & j & 0 & i & 0 & 0 & 0 & 0 & 0\\
 -j & 0 & i & 0 & 0 & 0 & 0 & 0 & 0\\
 0 & 0 & 0 & 0 & (1+k)/\sqrt{2} & 0 & 0 & 0 & 0\\
 0 & 0 & 0 & (1-k)/\sqrt{2} & 0 & (1+i)/\sqrt{2} & 0 & 0 & 0\\
 0 & 0 & 0 & 0 & 0 & 0 & -j & 0 & k\\
 0 & 0 & 0 & 0 & 0 & 1 & 0 & -k & 0\\
 \end{array}
 \right)\\
\\&=6.
 \end{align*}
\begin{align*}
 A(\widetilde{G}_{2})&= \left (
\begin{array}{ccccccccccccc}
 0 & j & 0 & 0 & 0 & -i & 0\\
 -j & 0 & k & 0 & 0 & 0 & 0\\
 0 & -k & 0 & k & 0 & 0 & 0\\
 0 & 0 & -k & 0 & j & 0 & -j\\
 0 & 0 & 0 & -j & 0 & i & 0\\
 i & 0 & 0 & 0 & -i & 0 & i\\
 0 & 0 & 0 & j & 0 & -i & 0\\
\end{array}
 \right).
 \end{align*}
We add 2-th row to 4-th row and add 5-th row to 7-th row:
\begin{align*}
 r(\widetilde{G}_{2})&=r \left(
\begin{array}{ccccccccccccc}
 0 & j & 0 & 0 & 0 & -i & 0\\
 -j & 0 & k & 0 & 0 & 0 & 0\\
 0 & -k & 0 & k & 0 & 0 & 0\\
 -j & 0 & 0 & 0 & j & 0 & -j\\
 0 & 0 & 0 & -j & 0 & i & 0\\
 i & 0 & 0 & 0 & -i & 0 & i\\
 0 & 0 & 0 & 0 & 0 & 0 & 0\\
 \end{array}
 \right).
 \end{align*}
We add 5-th row to 1-th row,
multiply $-k$ on the left side of 4-th row and add it to 6-th row:
\begin{align*}
  r(\widetilde{G}_{2})=r \left(
\begin{array}{ccccccccccccc}
 0 & j & 0 & -j & 0 & 0 & 0\\
 -j & 0 & k & 0 & 0 & 0 & 0\\
 0 & -k & 0 & k & 0 & 0 & 0\\
 -j & 0 & 0 & 0 & j & 0 & -j\\
 0 & 0 & 0 & -j & 0 & i & 0\\
 0 & 0 & 0 & 0 & 0 & 0 & 0\\
 0 & 0 & 0 & 0 & 0 & 0 & 0\\
  \end{array}
 \right).
 \end{align*}
We multiply $-i$ on the left side of 3-th row and add it to 1-th row,
\begin{align*}
  r(\widetilde{G}_{2})=r \left(
\begin{array}{ccccccccccccc}
 0 & 0 & 0 & 0 & 0 & 0 & 0\\
 -j & 0 & k & 0 & 0 & 0 & 0\\
 0 & -k & 0 & k & 0 & 0 & 0\\
 -j & 0 & 0 & 0 & j & 0 & -j\\
 0 & 0 & 0 & -j & 0 & i & 0\\
 0 & 0 & 0 & 0 & 0 & 0 & 0\\
 0 & 0 & 0 & 0 & 0 & 0 & 0\\
 \end{array}
 \right)
 =r \left(
\begin{array}{ccccccccccccc}
 -j & 0 & k & 0 & 0 & 0 & 0\\
 0 & -k & 0 & k & 0 & 0 & 0\\
 -j & 0 & 0 & 0 & j & 0 & -j\\
 0 & 0 & 0 & -j & 0 & i & 0\\
 \end{array}
 \right)=4.
 \end{align*}
Then $r(\widetilde{G}_{1})=6=|V(\widetilde{G}_{1})|-2c(\widetilde{G}_{1})+1$ and $r(\widetilde{G}_{2})=4=|V(\widetilde{G}_{2})|-2c(\widetilde{G}_{2})+1$.
Thus $\widetilde{G}_{1}$ and $\widetilde{G}_{2}$ are two extremal graphs which satisfy the condition in Theorem \ref{th:5.2.}.
\end{example}

\noindent\begin{lemma}\label{le:4.8.}
Let $\widetilde{G}=\widetilde{\theta}(p,l,q)$ be a $U(\mathbb{Q})$-gain graph of order $n$, where each cycle of $\widetilde{G}$ is of Type $1$ and $p,l,q$ are odd. Then $r(\widetilde{G}-v)=n-3$ for any vertex $v$ in $\widetilde{G}$.
\end{lemma}
\noindent\textbf{Proof.}
By Lemma \ref{le:2.4.} and Theorem \ref{th:5.2.}, $r(\widetilde{G}-v)\leq r(\widetilde{G})=n-3$.

\textbf{Case 1.} $d_{\widetilde{G}}(v)=2$.

If $\widetilde{G}-v$ is a $U(\mathbb{Q})$-gain cycle, by Lemma \ref{le:2.3.}, $r(\widetilde{G}-v)=|V(\widetilde{G}-v)|-2=n-3$.
If $\widetilde{G}-v$ has only one leaf and $c(\widetilde{G}-v)=1$, then the distance from the leaf to the $U(\mathbb{Q})$-gain cycle is even. By Lemma \ref{le:2.5.} several times, $r(\widetilde{G}-v)=n-3$.
Let $\widetilde{G}-v$ have two leaves and $c(\widetilde{G}-v)=1$.
If the distance from a leaf to the $U(\mathbb{Q})$-gain cycle is odd, by Lemma \ref{le:2.5.} and Theorem \ref{th:3.2.}, $r(\widetilde{G}-v)\geq n-3$, thus $r(\widetilde{G}-v)=n-3$.
If the distance from a leaf to the $U(\mathbb{Q})$-gain cycle is even, then, by Lemmas \ref{le:2.3.} and \ref{le:2.5.}, $r(\widetilde{G}-v)=n-3$.

\textbf{Case 2.} $d_{\widetilde{G}}(v)=3$.

Now, $\widetilde{G}-v$ is a $U(\mathbb{Q})$-gain tree, by Lemma \ref{le:2.5.}, $r(\widetilde{G}-v)\geq n-3$, thus $r(\widetilde{G}-v)=n-3$.
\quad $\square$~\\

\begin{figure}[htbp]
  \centering
  \includegraphics[scale=0.9]{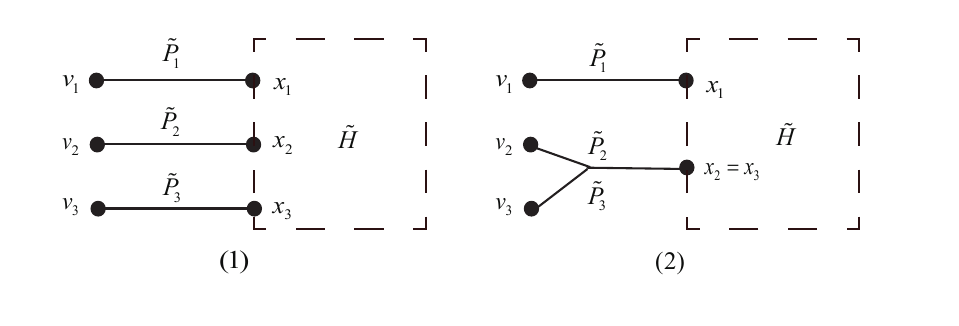}
  \caption{$\widetilde{G}-y$.}
\end{figure}

\noindent\begin{lemma}\label{le:4.9.}
Let $\widetilde{G}$ be a connected leaf-free $U(\mathbb{Q})$-gain graph of order $n$ and $c(\widetilde{G})\geq3$. If $r(\widetilde{G})=n-2c(\widetilde{G})+1$, then
\begin{enumerate}[(1)]
  \item $\widetilde{G}$ contains a cut-point and a pendant cycle;
  \item If $\widetilde{C}_{t}$ is a pendant cycle, $\widetilde{H}_{1}$ the maximal leaf-free subgraph of $\widetilde{G}-\widetilde{C}_{t}$ and $\widetilde{P}_{m}$ the path connecting $\widetilde{C}_{t}$ and $\widetilde{H}_{1}$, then $m$ is even.
\end{enumerate}
\end{lemma}
\noindent\textbf{Proof.}
Suppose that $\widetilde{G}$ does not contain any cut-point. Then for any vertex $v\in V(\widetilde{G})$, $\widetilde{G}-v$ is connected.
By Lemma \ref{le:3.1.} (3), we get $d_{\widetilde{G}}(v)=c(\widetilde{G})-c(\widetilde{G}-v)+1$.

If there is a vertex $x\in V(\widetilde{G})$ such that $d_{\widetilde{G}}(x)\geq4$, then $c(\widetilde{G}-x)\leq c(\widetilde{G})-3$.

When $\widetilde{G}-x$ does not contain any leaf, by Lemma \ref{le:2.4.} and Theorem \ref{th:3.2.}, we get
$$r(\widetilde{G})\geq r(\widetilde{G}-x)\geq |V(\widetilde{G}-x)|-2c(\widetilde{G}-x)\geq n-1-2(c(\widetilde{G})-3)= n-2c(\widetilde{G})+5,$$
a contradiction.

When $\widetilde{G}-x$ contains leaves, then $\widetilde{G}-x$ has at most $d_{\widetilde{G}}(x)$ leaves. By Lemma \ref{le:2.4.} and Theorem \ref{th:3.2.},
\begin{align}
r(\widetilde{G})&\geq r(\widetilde{G}-x)\nonumber
\\&\geq |V(\widetilde{G}-x)|-2c(\widetilde{G}-x)-p(\widetilde{G}-x)+1\nonumber
\\&\geq n-1-2c(\widetilde{G}-x)-d_{\widetilde{G}}(x)+1\nonumber
\\&=n-1-2c(\widetilde{G}-x)-(c(\widetilde{G})-c(\widetilde{G}-x)+1)+1\nonumber
\\&=n-c(\widetilde{G}-x)-c(\widetilde{G})-1\nonumber
\\&\geq n-(c(\widetilde{G})-3)-c(\widetilde{G})-1\nonumber
\\&=n-2c(\widetilde{G})+2,\nonumber
\end{align}
a contradiction. Then for any vertex $v\in V(\widetilde{G})$, $d_{\widetilde{G}}(v)\leq3$.

There must be a vertex $y\in V(\widetilde{G})$ such that $d_{\widetilde{G}}(y)=3$ (if all vertices of $\widetilde{G}$ with degree 2, then $\widetilde{G}$ is a $U(\mathbb{Q})$-gain cycle, which contradicts $c(\widetilde{G})\geq3$). Then $\widetilde{G}-y$ is a connected $U(\mathbb{Q})$-gain graph with $c(\widetilde{G}-y)=c(\widetilde{G})-2$.

When $\widetilde{G}-y$ does not contain any leaf, by Lemma \ref{le:2.4.} and Theorem \ref{th:3.2.},
$$r(\widetilde{G})\geq r(\widetilde{G}-y)\geq |V(\widetilde{G}-y)|-2c(\widetilde{G}-y)=(n-1)-2(c(\widetilde{G})-2)=n-2c(\widetilde{G})+3,$$
a contradiction.

When $\widetilde{G}-y$ contains leaves, then $\widetilde{G}-y$ has at most three leaves (since $\widetilde{G}$ does not contain any leaf). Now we discuss two cases.

\textbf{Case 1.} $\widetilde{G}-y$ contains at most two leaves.

By Lemma \ref{le:2.4.} and Theorem \ref{th:3.2.}, we get
\begin{align}
r(\widetilde{G})&\geq r(\widetilde{G}-y)\nonumber
\\&\geq |V(\widetilde{G}-y)|-2c(\widetilde{G}-y)-p(\widetilde{G}-y)+1\nonumber
\\&\geq(n-1)-2(c(\widetilde{G})-2)-2+1\nonumber
\\&=n-2c(\widetilde{G})+2,\nonumber
\end{align}
a contradiction.

\textbf{Case 2.} $\widetilde{G}-y$ contains three leaves $v_{1},v_{2}$ and $v_{3}$.

Let $\widetilde{H}$ be the maximal leaf-free subgraph of $\widetilde{G}-y$ (see Fig. 4) and let $\widetilde{P}_{i}$ $(i=1,2,3)$ be the unique $U(\mathbb{Q})$-gain path from $v_{i}$ to a vertex $x_{i}$ of $\widetilde{H}$.
Then $\widetilde{P}_{1}$, $\widetilde{P}_{2}$ and $\widetilde{P}_{3}$ do not contain a common vertex (since $\widetilde{G}$ does not contain any cut-point).
Let $\widetilde{P}_{1}$ be a $U(\mathbb{Q})$-gain path of order $m$ and vertex-disjoint from $\widetilde{P}_{2}$ and $\widetilde{P}_{3}$.

Since $\widetilde{G}$ does not contain any cut-point, if $\widetilde{G}-y$ is (1) of Fig. 4, then at least one of $x_{1}$, $x_{2}$ and $x_{3}$ lies on some cycle of $\widetilde{H}$. Without loss of generality, we assume that $x_{1}$ lies on a $U(\mathbb{Q})$-gain cycle of $\widetilde{H}$. Thus $c(\widetilde{H}-x_{1})\leq c(\widetilde{H})-1$. If $\widetilde{G}-y$ is (2) of Fig. 4, then each $x_{i}$ $(i=1,2,3)$ lies on a $U(\mathbb{Q})$-gain cycle of $\widetilde{H}$. Thus $c(\widetilde{H}-x_{1})\leq c(\widetilde{H})-1$.

When $m$ is even, let $\widetilde{G}_{1}=\widetilde{G}-y-\widetilde{P}_{1}$. Then $\widetilde{G}_{1}$ contains at most four leaves and $c(\widetilde{G}_{1})=c(\widetilde{H}-x_{1})\leq c(\widetilde{H})-1$. By Lemmas \ref{le:2.4.}, \ref{le:2.5.} and Theorem \ref{th:3.2.}, we obtain
\begin{align}
r(\widetilde{G})&\geq r(\widetilde{G}-y)=r(\widetilde{G}_{1})+m\nonumber
\\&\geq |V(\widetilde{G}_{1})|-2c(\widetilde{G}_{1})-p(\widetilde{G}_{1})+1+m\nonumber
\\&\geq (n-m-1)-2(c(\widetilde{H})-1)-4+1+m\nonumber
\\&=n-2c(\widetilde{H})-2=n-2c(\widetilde{G}-y)-2\nonumber
\\&= n-2(c(\widetilde{G})-2)-2=n-2c(\widetilde{G})+2,\nonumber
\end{align}
a contradiction.

When $m$ is odd, let $\widetilde{G}_{2}=\widetilde{G}-y-\widetilde{P}_{1}+x_{1}$. Then $\widetilde{G}_{2}$ contains two leaves and $c(\widetilde{G}_{2})=c(\widetilde{G}-y)$. By Lemmas \ref{le:2.4.}, \ref{le:2.5.} and Theorem \ref{th:3.2.},
\begin{align}
r(\widetilde{G})&\geq r(\widetilde{G}-y)=r(\widetilde{G}_{2})+m-1\nonumber
\\&\geq |V(\widetilde{G}_{2})|-2c(\widetilde{G}_{2})-p(\widetilde{G}_{2})+1+m-1\nonumber
\\&=n-2c(\widetilde{G}-y)-2= n-2(c(\widetilde{G})-2)-2=n-2c(\widetilde{G})+2,\nonumber
\end{align}
a contradiction.

Thus $\widetilde{G}$ contains a cut-point; next we will prove that $\widetilde{G}$ contains a pendant cycle.

Let $u$ be a cut-point of $\widetilde{G}$ and let $\widetilde{G}$ be obtained from a block $\widetilde{G}_{3}$ and a $U(\mathbb{Q})$-gain graph $\widetilde{G}_{4}$ by identifying the unique common vertex $u$. Since $\widetilde{G}$ is a connected leaf-free $U(\mathbb{Q})$-gain graph and $\widetilde{G}_{3}$ is a leaf-free block of $\widetilde{G}$, $\widetilde{G}_{4}$ has at most one leaf. By Theorem \ref{th:3.2.}, $r(\widetilde{G}_{4})\geq|V(\widetilde{G}_{4})|-2c(\widetilde{G}_{4})$.

If $c(\widetilde{G}_{3})=1$, then, by Lemma \ref{le:4.5.}, either $\widetilde{G}_{3}$ is a pendant cycle of $\widetilde{G}$, or $\widetilde{G}=\widetilde{\infty}(p,1,q)$, where each cycle of $\widetilde{G}$ is of Type $1$, a contradiction.

If $c(\widetilde{G}_{3})=2$, then $\widetilde{G}_{3}=\widetilde{\theta}(p,l,q)$. By Theorem \ref{th:3.2.}, $r(\widetilde{G}_{3})\geq|V(\widetilde{G}_{3})|-2c(\widetilde{G}_{3})+1$.
If $r(\widetilde{G}_{3})\geq |V(\widetilde{G}_{3})|-2c(\widetilde{G}_{3})+2$, then, by Lemma \ref{le:4.3.}, $r(\widetilde{G})\geq |V(\widetilde{G})|-2c(\widetilde{G})+2$, a contradiction.
Thus $r(\widetilde{G}_{3})=|V(\widetilde{G}_{3})|-2c(\widetilde{G}_{3})+1$; by Theorem \ref{th:5.2.}, each $U(\mathbb{Q})$-gain cycle of $\widetilde{G}_{3}$ is of Type 1 and $p,l,q$ are odd. Combining with Lemmas \ref{le:2.10.}, \ref{le:2.4.} and \ref{le:2.5.}, we get
\begin{align}
r(\widetilde{G})&\geq p-1+r(\widetilde{C}_{l+q+2})+r(\widetilde{G}_{4})=|V(\widetilde{G}_{3})|-3+r(\widetilde{G}_{4})\nonumber
\\&\geq |V(\widetilde{G}_{3})|-3+(|V(\widetilde{G}_{4})|-2c(\widetilde{G}_{4}))=n-2c(\widetilde{G})+2,\nonumber
\end{align}
a contradiction.

If $c(\widetilde{G}_{3})\geq3$, then, by Theorem \ref{th:3.2.}, $r(\widetilde{G}_{3})\geq|V(\widetilde{G}_{3})|-2c(\widetilde{G}_{3})+1$.
If $r(\widetilde{G}_{3})=|V(\widetilde{G}_{3})|-2c(\widetilde{G}_{3})+1$, then $\widetilde{G}_{3}$ contains a cut-point, a contradiction.
If $r(\widetilde{G}_{3})\geq |V(\widetilde{G}_{3})|-2c(\widetilde{G}_{3})+2$, then, by Lemma \ref{le:4.3.}, $r(\widetilde{G})\geq n-2c(\widetilde{G})+2$, a contradiction.

According to the above discussion, (1) holds.

By Lemma \ref{le:4.7.}, $\widetilde{C}_{t}$ is of Type 1. Let $v$ be a major vertex which lies on $\widetilde{C}_{t}$. If $m$ is odd, then, by Lemmas \ref{le:2.5.} and \ref{le:2.10.},
\begin{align}
r(\widetilde{G})&=t-2+r(\widetilde{G}-\widetilde{C}_{t}+v)\nonumber
\\&=t-2+m-1+r(\widetilde{H}_{1})\nonumber
\\&\geq t-2+m-1+(|V(\widetilde{H}_{1})|-2c(\widetilde{H}_{1})+1)\nonumber
\\&=t-2+m-1+|V(\widetilde{H}_{1})|-2(c(\widetilde{G})-1)+1\nonumber
\\&=n-2c(\widetilde{G})+2,\nonumber
\end{align}
a contradiction. Then $m$ is even, hence (2) holds.
\quad $\square$~\\

Now, we can prove the main result of this section.

\begin{figure}[!ht]
  \centering
 \includegraphics[scale=0.7]{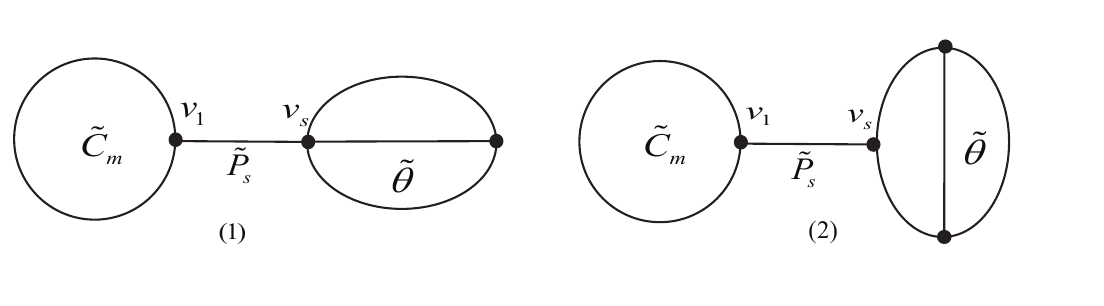}
 \caption{$\widetilde{G}$.}
\end{figure}
\begin{figure}[htbp]
\centering
 \includegraphics[scale=0.7]{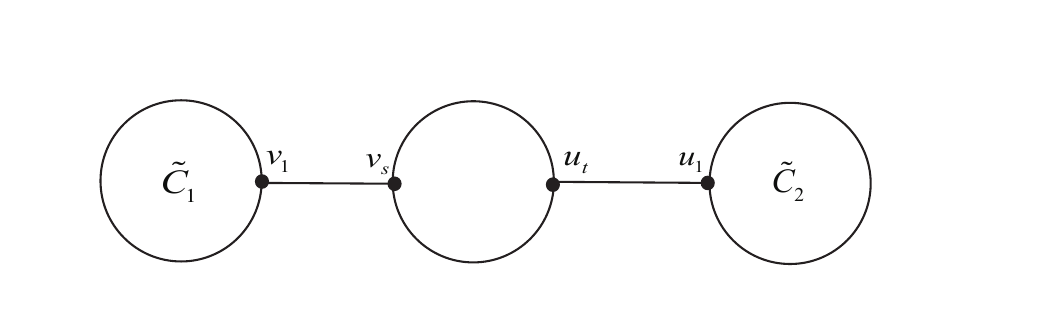}
 \caption{$\widetilde{G}$.}
\end{figure}

\noindent\begin{theorem}\label{th:5.3.}
Let $\widetilde{G}$ be a $U(\mathbb{Q})$-gain graph of order $n$ with $p(\widetilde{G})=0$ and $c(\widetilde{G})\geqslant3$. Then $r(\widetilde{G})=n-2c(\widetilde{G})+1$ if and only if $\widetilde{G}$ is a $U(\mathbb{Q})$-gain graph which is obtained from a tree $\widetilde{T}$ with $p(\widetilde{T})=c(\widetilde{G})$ and $r(\widetilde{T})=|V(\widetilde{T})|-p(\widetilde{T})+1$, by attaching a cycle which is of Type $1$ on each leaf of $\widetilde{T}$.
\end{theorem}
\noindent\textbf{Proof.}
\textbf{Sufficiency:}
Since $\widetilde{G}$ is obtained from a tree $\widetilde{T}$ with $p(\widetilde{T})=c(\widetilde{G})=c$ and $r(\widetilde{T})=|V(\widetilde{T})|-p(\widetilde{T})+1$, by attaching a cycle which is of Type $1$ on each leaf of $\widetilde{T}$, then each cycle becomes a pendant cycle of $\widetilde{G}$.
Let $\widetilde{C}_{1}, \widetilde{C}_{2}, \ldots, \widetilde{C}_{c}$ be pendant cycles of $\widetilde{G}$. By Lemma \ref{le:2.10.}, we get
$$r(\widetilde{G})=\sum^{c}_{i=1}\{|V(\widetilde{C}_{i})|-2\}+r(\widetilde{T})=\sum^{c}_{i=1}|V(\widetilde{C}_{i})|-2c+|V(\widetilde{T})|-c+1=n-2c+1.$$

\textbf{Necessity:}
Since $r(\widetilde{G})=n-2c(\widetilde{G})+1$, by Lemmas \ref{le:4.7.} and \ref{le:4.9.}$(1)$, $\widetilde{G}$ contains a pendant cycle which is of Type $1$. Next, we will apply induction on $c(\widetilde{G})$.

When $c(\widetilde{G})=3$, then $r(\widetilde{G})=n-5$. We first prove that all cycles of $\widetilde{G}$ are pendant cycles.

If $\widetilde{G}$ has only one pendant cycle $\widetilde{C}_{m}$, then $\widetilde{G}$ is obtained from $\widetilde{C}_{m}$ and $\widetilde{\theta}=\widetilde{\theta}(p,l,q)$ by connecting $v_{1}\in V(\widetilde{C}_{m})$ and $v_{s}\in V(\widetilde{\theta})$ with a path $\widetilde{P}(v_{1},v_{2},\ldots,v_{s})$, $s\geq2$ (see Fig. 5). By Lemma \ref{le:4.9.}(2), $s$ is even, then $r(\widetilde{G}-\widetilde{C}_{m}+v_{1})=r(\widetilde{\theta}-v_{s})+s$ by Lemma \ref{le:2.5.}. Combining with Lemma \ref{le:2.10.}, we get
$$r(\widetilde{G})=m-2+r(\widetilde{G}-\widetilde{C}_{m}+v_{1})=m-2+s+r(\widetilde{\theta}-v_{s}).$$

When $\widetilde{G}$ is (1) in Fig. 5, $\widetilde{\theta}-v_{s}$ is a tree with $1\leq p(\widetilde{\theta}-v_{s})\leq3$. By Theorem \ref{th:3.2.}, $r(\widetilde{\theta}-v_{s})\geq |V(\widetilde{\theta}-v_{s})|-2$. Then $r(\widetilde{G})\geq n-3$, a contradiction.

When $\widetilde{G}$ is (2) in Fig. 5,
$\widetilde{\theta}-v_{s}$ is a connected $U(\mathbb{Q})$-gain graph with $c(\widetilde{\theta}-v_{s})=1$ and $1\leq p(\widetilde{\theta}-v_{s})\leq2$, then, by Theorem \ref{th:3.2.}, $r(\widetilde{\theta}-v_{s})\geq |V(\widetilde{\theta}-v_{s})|-3$. Thus $r(\widetilde{G})\geq n-4$, a contradiction.
If $\widetilde{\theta}-v_{s}$ is a $U(\mathbb{Q})$-gain cycle, then, by Theorem \ref{th:3.2.}, $r(\widetilde{\theta}-v_{s})\geq |V(\widetilde{\theta}-v_{s})|-2$. Thus $r(\widetilde{G})\geq n-3$, a contradiction.

If $\widetilde{G}$ has two pendant cycles $\widetilde{C}_{1}$ and $\widetilde{C}_{2}$, then $\widetilde{G}$ is as shown in Fig. 6.

By Lemma \ref{le:4.9.}(2), $s$ and $t$ are even. Let $\widetilde{G}_{1}=\widetilde{G}-\widetilde{C}_{1}-\widetilde{C}_{2}+v_{1}+u_{1}$. By Lemma \ref{le:2.5.} and Theorem \ref{th:3.2.}, $r(\widetilde{G}_{1})\geq |V(\widetilde{G}_{1})|-2$. Combining with Lemma \ref{le:2.10.}, we obtain
$$r(\widetilde{G})=|V(\widetilde{C}_{1})|+|V(\widetilde{C}_{2})|-4+r(\widetilde{G}_{1})\geq n-4,$$
a contradiction. So, all cycles of $\widetilde{G}$ are pendant cycles.

Let $\widetilde{C}_{1}$, $\widetilde{C}_{2}$ and $\widetilde{C}_{3}$ be all pendant cycles of $\widetilde{G}$. Let $\widetilde{T}$ be a $U(\mathbb{Q})$-gain tree obtained by shrinking $\widetilde{C}_{i}(i=1,2,3)$ into a vertex. By Lemma \ref{le:2.10.}, we obtain
$$r(\widetilde{G})=|V(\widetilde{C}_{1})|+|V(\widetilde{C}_{2})|+|V(\widetilde{C}_{3})|-6+r(\widetilde{T}).$$
Since $r(\widetilde{G})=n-5$, $r(\widetilde{T})=|V(\widetilde{T})|-2=|V(\widetilde{T})|-p(\widetilde{T})+1$. Thus, $\widetilde{G}$ is obtained from a tree $\widetilde{T}$ with $p(\widetilde{T})=c(\widetilde{G})=3$ and $r(\widetilde{T})=|V(\widetilde{T})|-p(\widetilde{T})+1$, by attaching a cycle which is of Type $1$ on each leaf of $\widetilde{T}$.

We assume that the result holds when $c(\widetilde{G})<c$. Now we consider the case $c(\widetilde{G})=c (c\geq4)$.

Since $r(\widetilde{G})=n-2c(\widetilde{G})+1$, $\widetilde{G}$ has a pendant cycle $\widetilde{C}_{m}$ which is of Type $1$ (by Lemma \ref{le:4.7.} and \ref{le:4.9.}(1)). Let $\widetilde{H}$ be a maximal leaf-free subgraph of $\widetilde{G}-\widetilde{C}_{m}$, then $\widetilde{G}$ is obtained from $\widetilde{C}_{m}$ and $\widetilde{H}$ by connecting $v_{1}\in V(\widetilde{C}_{m})$ and $v_{s}\in V(\widetilde{H})$ with a path $\widetilde{P}(v_{1},v_{2},\ldots,v_{s}) (s\geqslant2)$ and $s$ is even (by Lemma \ref{le:4.9.}(2)).

Let $\widetilde{G}_{1}=\widetilde{G}-\widetilde{C}_{m}+v_{1}$. By Lemma \ref{le:2.10.}, we get
$$r(\widetilde{G})=m-2+r(\widetilde{G}_{1}).$$
Since $r(\widetilde{G})=n-2c(\widetilde{G})+1$, $r(\widetilde{G}_{1})=|V(\widetilde{G}_{1})|-2c(\widetilde{G}_{1})$. By Lemmas \ref{le:2.4.} and \ref{le:2.5.},
$$r(\widetilde{H})\leq r(\widetilde{H}+v_{s-1})=r(\widetilde{G}_{1})-(s-2)=|V(\widetilde{H})|-2c(\widetilde{H})+1.$$
Since $c(\widetilde{G})\geq3$, $\widetilde{H}$ is not a cycle. By Theorems \ref{th:3.2.} and \ref{th:5.1.}, $r(\widetilde{H})\geq |V(\widetilde{H})|-2c(\widetilde{H})+1$, thus $r(\widetilde{H})=|V(\widetilde{H})|-2c(\widetilde{H})+1$. Since $c(\widetilde{H})=c(\widetilde{G})-1<c$, by the induction assumption, $\widetilde{H}$ is obtained from a tree $\widetilde{T}_{1}$ with $p(\widetilde{T}_{1})=c-1$ and $r(\widetilde{T}_{1})=|V(\widetilde{T}_{1})|-p(\widetilde{T}_{1})+1$, by attaching a cycle which is of Type $1$ on each leaf of $\widetilde{T}_{1}$. So all cycles of $\widetilde{H}$ are pendant cycles, thus $\widetilde{G}$ has at least $c-1$ pendant cycles.

If $\widetilde{G}$ has $c-1$ pendant cycles, then $v_{s}$ is on a $U(\mathbb{Q})$-gain cycle of $\widetilde{H}$, $c(\widetilde{H}-v_{s})=c(\widetilde{G}_{1})-1$ and $p(\widetilde{H}-v_{s})\leq 2$. By Lemma \ref{le:2.5.}, we have
$$r(\widetilde{H}-v_{s})=r(\widetilde{G}_{1})-s=|V(\widetilde{H}-v_{s})|-2c(\widetilde{H}-v_{s})-2.$$
By Theorem \ref{th:3.2.}, $r(\widetilde{H}-v_{s})\geq |V(\widetilde{H}-v_{s})|-2c(\widetilde{H}-v_{s})-1$, a contradiction.

Therefore, $\widetilde{G}$ has $c$ pendant cycles, and each cycle is of Type $1$. Let $\widetilde{C}_{1}, \widetilde{C}_{2}, \ldots, \widetilde{C}_{c}$ be pendant cycles of $\widetilde{G}$. By Lemma \ref{le:2.10.}, we have
$$r(\widetilde{G})=\sum^{c}_{i=1}\{|V(\widetilde{C}_{i})|-2\}+r(\widetilde{T})=\sum^{c}_{i=1}|V(\widetilde{C}_{i})|-2c+r(\widetilde{T}).$$
Since $r(\widetilde{G})=n-2c(\widetilde{G})+1$, $r(\widetilde{T})=|V(\widetilde{T})|-c+1$. Thus, $\widetilde{G}$ is obtained from a tree $\widetilde{T}$ with $p(\widetilde{T})=c(\widetilde{G})$ and $r(\widetilde{T})=|V(\widetilde{T})|-p(\widetilde{T})+1$, by attaching a cycle which is of Type $1$ on each leaf of $\widetilde{T}$.
\quad $\square$

\begin{figure}[htbp]
  \centering
 \includegraphics[scale=0.7]{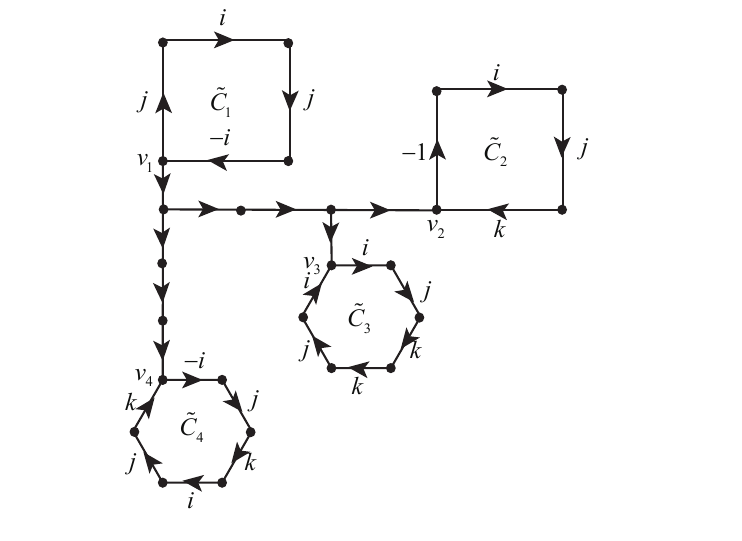}
  \caption{$\widetilde{G}$.}
\end{figure}

\noindent\begin{example}
Consider the $U(\mathbb{Q})$-gain graph $\widetilde{G}$ in Fig. 7, then
$$\varphi(\widetilde{C}_{1})=ij(-i)j=(ij)(-ij)=k(-k)=1;$$
$$\varphi(\widetilde{C}_{2})=(-1)ijk=-k^{2}=1;$$
$$\varphi(\widetilde{C}_{3})=ijkkji=(ij)(kk)(ji)=k(-1)(-k)=-1;$$
$$\varphi(\widetilde{C}_{4})=-ijkijk=-(ij)k(ij)k=-k^{4}=-1.$$
By Definition \ref{de:2.2.}, $\widetilde{C}_{1}, \widetilde{C}_{2}, \widetilde{C}_{3}$ and $\widetilde{C}_{4}$ are of Type 1.
Let $\widetilde{G}'=\widetilde{G}-\widetilde{C}_{1}-\widetilde{C}_{2}-\widetilde{C}_{3}-\widetilde{C}_{4}+v_{1}+v_{2}+v_{3}+v_{4}$, then $\widetilde{G}'$ is a $U(\mathbb{Q})$-gain tree $\widetilde{T}$. By Lemma \ref{le:2.5.}, we get $r(\widetilde{T})=6=|V(\widetilde{T})|-p(\widetilde{T})+1$.
Combining with Lemma \ref{le:2.10.}, we obtain
$$r(\widetilde{G})=(4-2)+(4-2)+(6-2)+(6-2)+r(\widetilde{T})=18=|V(\widetilde{G})|-2c(\widetilde{G})+1.$$
Thus $\widetilde{G}$ is an extremal graph which satisfies the condition in Theorem \ref{th:5.3.}.
\end{example}

\section{$U(\mathbb{Q})$-gain graphs $\widetilde{G}$ with $r(\widetilde{G})=|V(\widetilde{G})|-2c(\widetilde{G})-p(\widetilde{G})+1$}

Lastly, we give a characterization of the connected $U(\mathbb{Q})$-gain graphs $\widetilde{G}$ with $p(\widetilde{G})\geqslant1$ and $r(\widetilde{G})=|V(\widetilde{G})|-2c(\widetilde{G})-p(\widetilde{G})+1$.
When $c(\widetilde{G})=0$ and $p(\widetilde{G})=2$, then $\widetilde{G}=\widetilde{P}_{t}$ and by Lemma \ref{le:2.1.}, the equality holds if and only if $t$ is odd.
When $c(\widetilde{G})=0$ and $p(\widetilde{G})\geq3$, then $\widetilde{G}=\widetilde{T}$ is a tree.
In this situation, the graph with $r(\widetilde{T})=|V(\widetilde{T})|-p(\widetilde{T})+1$ is characterized in Theorem \ref{th:5.4.}.
Next, we give a characterization of $U(\mathbb{Q})$-gain graphs with $c(\widetilde{G})\geq1$, $p(\widetilde{G})\geqslant1$ and $r(\widetilde{G})=|V(\widetilde{G})|-2c(\widetilde{G})-p(\widetilde{G})+1$.

\noindent\begin{theorem}\label{th:5.5.}
Let $\widetilde{G}$ be a connected $U(\mathbb{Q})$-gain graph of order $n$ with $c(\widetilde{G})\geqslant1$ and $p(\widetilde{G})\geqslant1$. Then $r(\widetilde{G})=n-2c(\widetilde{G})-p(\widetilde{G})+1$ if and only if $\widetilde{G}$ is obtained form a tree $\widetilde{T}$ with $p(\widetilde{T})>c(\widetilde{G})$ and $r(\widetilde{T})=|V(\widetilde{T})|-p(\widetilde{T})+1$, by attaching $c(\widetilde{G})$ cycles on $c(\widetilde{G})$ leaves of $\widetilde{T}$, each cycle of Type $1$.
\end{theorem}
\noindent\textbf{Proof.}
\textbf{Sufficiency:}
Let $\widetilde{C}_{1}, \widetilde{C}_{2}, \ldots, \widetilde{C}_{c}$ be all pendant cycles of $\widetilde{G}$. By Lemma \ref{le:2.10.}, we get
\begin{align}
r(\widetilde{G})&=\sum^{c}_{i=1}\{|V(\widetilde{C}_{i})|-2\}+r(\widetilde{T})=\sum^{c}_{i=1}|V(\widetilde{C}_{i})|-2c+|V(\widetilde{T})|-p(\widetilde{T})+1\nonumber
\\&=n-2c(\widetilde{G})-p(\widetilde{G})+1.\nonumber
\end{align}

\textbf{Necessity:}
Let $B$ be a block of $\widetilde{G}$ such that $c(B)\geq c(B')$ for any block $B'$ of $\widetilde{G}$. Then $r(B)=|V(B)|-2c(B)-p(B)+1$ (by Lemma \ref{le:4.3.}). Thus by Lemma \ref{le:4.9.}, $c(B)=1$ or $2$.

\textbf{Case 1.} $c(B)=2$

Since $B$ is a block, $B=\widetilde{\theta}(p,l,q)$. By Theorem \ref{th:5.2.}, each cycle of $B$ is of Type $1$, and $p, l, q$ are odd. Let $x$ be a cut-point of $\widetilde{G}$ which lies on $B$. By Lemma \ref{le:4.8.}, $r(B-x)=|V(B)|-3$. Let $\widetilde{H}_{1}$ be a component of $\widetilde{G}-x$ such that $V(\widetilde{H}_{1})\cap V(B)=\emptyset$, $\widetilde{H}=\widetilde{H}_{1}+x$ and let $\widetilde{G}_{1}=B+\widetilde{H}_{1}$. Combining with Lemmas \ref{le:2.10.}, \ref{le:2.4.} and \ref{le:2.5.}, we obtain
$$r(\widetilde{G}_{1})\geq p-1+r(\widetilde{C}_{l+q+2})+r(\widetilde{H})=|V(B)|-3+r(\widetilde{H}).$$

If $\widetilde{H}$ is a $U(\mathbb{Q})$-gain cycle $\widetilde{C}_{t}$, then, by Lemma \ref{le:2.3.}, $r(\widetilde{H})\geq |V(\widetilde{H})|-2$. So $r(\widetilde{G}_{1})\geq|V(\widetilde{G}_{1})|-2c(\widetilde{G}_{1})-p(\widetilde{G}_{1})+2$.
If $\widetilde{H}$ is not a $U(\mathbb{Q})$-gain cycle, then, by Lemma \ref{le:2.3.}, $r(\widetilde{H})\geq |V(\widetilde{H})|-2c(\widetilde{H})-p(\widetilde{H})+1$. Note that $c(\widetilde{H})=c(\widetilde{G}_{1})-2$ and $p(\widetilde{H})\leq p(\widetilde{G}_{1})+1$, then
$$r(\widetilde{G}_{1})\geq|V(\widetilde{G}_{1})|-2c(\widetilde{G}_{1})-p(\widetilde{G}_{1})+2.$$

By Lemma \ref{le:4.3.}, $r(\widetilde{G})\geq n-2c(\widetilde{G})-p(\widetilde{G})+2$, a contradiction.

\textbf{Case 2.} $c(B)=1$

Let $\widetilde{C}_{1},\widetilde{C}_{2},\ldots, \widetilde{C}_{c}$ be all cycles of $\widetilde{G}$ and let $\widetilde{T}$ be a tree obtained by shrinking $\widetilde{C}_{i}~(i=1,2,\ldots,c)$ into a vertex. In this case, $\widetilde{C}_{i}~(i=1,2,\ldots,c)$ is a block. By Lemmas \ref{le:4.5.}, \ref{le:4.6.} and $p(\widetilde{G})\geqslant1$, we have that $\widetilde{C}_{i}$ is a pendant cycle of $\widetilde{G}$ and it is of Type $1$.

By Lemma \ref{le:2.10.}, we get $r(\widetilde{G})=\sum^{c}_{i=1}\{|V(\widetilde{C}_{i})|-2\}+r(\widetilde{T})$. Since $r(\widetilde{G})=n-2c-p(\widetilde{G})+1$, $r(\widetilde{T})=|V(\widetilde{T})|-p(\widetilde{T})+1$.
Thus $\widetilde{G}$ is obtained from a $U(\mathbb{Q})$-gain tree $\widetilde{T}$ with $p(\widetilde{T})>c(\widetilde{G})$ and $r(\widetilde{T})=|V(\widetilde{T})|-p(\widetilde{T})+1$, by attaching $c$ cycles which are of Type $1$, on $c$ leaves of $\widetilde{T}$.
\quad $\square$~\\

\begin{figure}[htbp]
  \centering
  \includegraphics[scale=0.75]{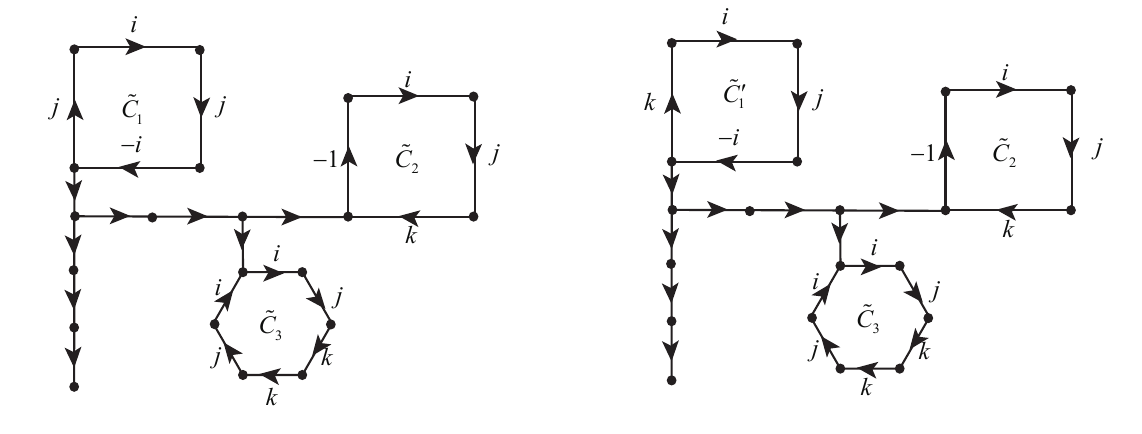}
  \caption{$\widetilde{G}_{1}$ and $\widetilde{G}_{2}$.}
\end{figure}

\noindent\begin{example}
Let $\widetilde{G}_{1}$ (see Fig. 8) be the $U(\mathbb{Q})$-gain graph obtained by shrinking $\widetilde{C}_{4}$ (in Example 5.12) into a vertex.
By Lemma \ref{le:2.10.}, we obtain
$$r(\widetilde{G}_{1})=(4-2)+(4-2)+(6-2)+r(\widetilde{T})=14=|V(\widetilde{G}_{1})|-2c(\widetilde{G}_{1})-p(\widetilde{G}_{1})+1.$$
Thus $\widetilde{G}_{1}$ is an extremal graph which satisfies the condition in Theorem \ref{th:5.5.}.

Let $\widetilde{G}_{2}$ (see Fig. 8) be the $U(\mathbb{Q})$-gain graph obtained by replacing $\widetilde{C}_{1}$ with $\widetilde{C}'_{1}$, where $\varphi(\widetilde{C}'_{1})=ij(-i)k=-i$.
By Definition \ref{de:2.2.}, $\widetilde{C}'_{1}$ is of Type 2.
Let $\widetilde{G}'_{2}=\widetilde{T}-v_{1}$.
Combining with Lemma \ref{le:2.10.}, we obtain
$$r(\widetilde{G}_{2})=4+(4-2)+(6-2)+r(\widetilde{G}'_{2})=16>|V(\widetilde{G}_{2})|-2c(\widetilde{G}_{2})-p(\widetilde{G}_{2})+1.$$
Thus $\widetilde{G}_{2}$ is a $U(\mathbb{Q})$-gain graph which does not satisfy the condition in Theorem \ref{th:5.5.}.
\end{example}

\end{document}